\documentclass[12pt]{amsart}
\usepackage{amsmath}
\usepackage{amsthm}
\usepackage{amssymb}
\newtheorem{theorem}{Theorem}[section]
\newtheorem{lemma}[theorem]{Lemma}
\newtheorem{proposition}[theorem]{Proposition}

\newtheorem{corollary}[theorem]{Corollary}
\newtheorem{remar}[theorem]{Remark}
\newtheorem{note}[theorem]{Notes}
\newtheorem{prob}{Open Problem}

\newenvironment{remark}{\begin{remar}\rm}{\end{remar}}

\newenvironment{axiom}{\begin{list}{$\bullet$}{\setlength{\labelsep}{.7cm}%
\setlength{\leftmargin}{2.5cm}\setlength{\rightmargin}{0cm}%
\setlength{\labelwidth}{1.8cm}\setlength{\itemsep}{0pt}}}{\end{list}}
\newcommand{\ax}[1]{\item[{\bf #1}\hfill]\index{#1}}
\newcommand{\bfind}[1]{\index{#1}{\bf #1}}

\newcommand{\n}{\par\noindent}

\newcommand{\sn}{\par\smallskip\noindent}
\newcommand{\mn}{\par\medskip\noindent}

\newcommand{\pars}{\par\smallskip}
\newcommand{\parm}{\par\medskip}

\newcommand{\isom}{\simeq}

\newcommand{\ind}[1]{\index{#1}{#1}}

\newcommand{\ec}{\prec_{\exists}}
\newcommand{\sep}{^{\rm sep}}

\newcommand{\chara}{\mbox{\rm char}\,}
\newcommand{\trdeg}{\mbox{\rm trdeg}\,}

\newcommand{\cal}{\mathcal}
\font\teneu=eufm10 scaled 1200
\font\seveneu=eufm7 scaled 1200
\font\fiveeu=eufm5 scaled 1200
\font\tenlv=msbm10 scaled 1200
\font\sevenlv=msbm7 scaled 1200
\font\fivelv=msbm5 scaled 1200

\def\eu #1{{\mathchoice{{\hbox{\teneu #1}}}{{\hbox{\teneu #1}}}
{{\hbox{\seveneu #1}}}{{\hbox{\fiveeu #1}}}}}
\def\lv #1{{\mathchoice{{\hbox{\tenlv #1}}}{{\hbox{\tenlv #1}}}
{{\hbox{\sevenlv #1}}}{{\hbox{\fivelv #1}}}}}

\newcommand{\N}{\lv N}
\newcommand{\Q}{\lv Q}

\newcommand{\Z}{\lv Z}
\newcommand{\F}{\lv F}

\begin{document}
\title{The algebra and model theory of tame valued fields}
\author{Franz--Viktor Kuhlmann}
\date{14.\ 3.\ 2014}
\thanks{The author would like to thank Peter Roquette, Alexander Prestel
and Florian Pop for very helpful discussions, and Koushik Pal for his
suggestions for improvement of an earlier version of the manuscript.
Special thanks to Anna Blaszczok for her amazingly careful and efficient
proofreading of the final version.\n
This research was partially supported by a Canadian NSERC grant.\n
AMS Subject Classification: 12J10, 12J15}
\address{Department of Mathematics and Statistics,
University of Saskatchewan, 106 Wiggins Road,
Saskatoon, Saskatchewan, Canada S7N 5E6}
\email{fvk@math.usask.ca}

\begin{abstract}\noindent
A henselian valued field $K$ is called a tame field if its algebraic
closure $\tilde{K}$ is a tame extension, that is, the ramification field
of the normal extension $\tilde{K}|K$ is algebraically closed. Every
algebraically maximal Kaplansky field is a tame field, but not
conversely. We develop the algebraic theory of tame fields and then
prove Ax--Kochen--Ershov Principles for tame fields. This
leads to model completeness and completeness results
relative to value group and residue field. As the maximal immediate
extensions of tame fields will in general not be unique, the proofs
have to use much deeper valuation theoretical results than those
for other classes of valued fields which have already been shown to
satisfy Ax--Kochen--Ershov Principles. The results of this paper have
been applied to gain insight in the Zariski space of places of an
algebraic function field, and in the model theory of large fields.
\end{abstract}

\maketitle

%
%
%
%
\section{Introduction}
In this paper, we consider valued fields. By $(K,v)$ we mean a field $K$
equipped with a valuation $v$. We denote the value group by $vK$, the
residue field by $Kv$ and the valuation ring by
${\cal O}_K\,$. For elements $a\in K$, the value is denoted by $va$, and
the residue by $av$. When we talk of a valued field extension $(L|K,v)$
we mean that $(L,v)$ is a valued field, $L|K$ a field extension, and $K$
is endowed with the restriction of $v$.
%
%

We write a valuation in the classical additive (Krull) way, that is, the
value group is an additively written ordered abelian group, the
homomorphism property of $v$ says that $vab= va+vb$, and the ultrametric
triangle law says that $v(a+b)\geq\min \{va,vb\}$. Further, we have the
rule $va=\infty \Leftrightarrow a=0$.

For the basic facts from valuation theory, we refer the reader to
\cite{[En]}, \cite{[E--P]}, \cite{[Ri]}, \cite{[W]}, \cite{[Z--S]} and
\cite{[K2]}.

\pars
In this paper, our main concern is the algebra and the model theory of
tame and of separably tame valued fields, which we will introduce now.

A valued field is \bfind{henselian} if it satisfies Hensel's Lemma, or
equivalently, if it admits a unique extension of the valuation to every
algebraic extension field. Take a henselian field
$(K,v)$, and let $p$ denote the \bfind{characteristic exponent} of its
residue field $Kv$, i.e., $p=\chara Kv$ if this is positive, and $p=1$
otherwise. An algebraic extension $(L|K,v)$ of a henselian field
$(K,v)$ is called \bfind{tame} if every finite subextension $E|K$ of
$L|K$ satisfies the following
conditions:
\begin{axiom}
\ax{(TE1)} The ramification index $(vE:vK)$ is prime to $p$,
\ax{(TE2)} The residue field extension $Ev|Kv$ is separable,
\ax{(TE3)} The extension $(E|K,v)$ is \bfind{defectless}, i.e.,
$[E:K]=(vE:vK)[Ev:Kv]$.
\end{axiom}
\begin{remark}
This notion of ``tame extension'' does not coincide with the notion of
``tamely ramified extension'' as defined on page 180 of O.~Endler's book
\cite{[En]}.
The latter definition requires (TE1) and
(TE2), but not (TE3). Our tame extensions are the defectless tamely
ramified extensions in the sense of Endler's book. In particular, in our
terminology, proper immediate algebraic extensions of henselian fields
are not called tame (in fact, they cause a lot of problems in the model
theory of valued fields).
\end{remark}

A \bfind{tame valued field} (in short, \bfind{tame field}) is a
henselian field for which all algebraic extensions are tame.
%
%
Likewise, a \bfind{separably tame field} is a henselian field
for which all separable-algebraic extensions are tame.
%
%
The algebraic properties and characterizations of tame fields will be
studied in Section~\ref{sectatf}, and those of separably tame fields in
Section~\ref{sectastf}.

If $\chara Kv=0$, then conditions (TE1) and (TE2) are void, and
every finite extension of $(K,v)$ is defectless (cf.\
Corollary~\ref{lostr0} below). We obtain the following well known fact:

\begin{theorem}                               \label{tamerc0}
Every algebraic extension of a henselian field of residue characteristic
0 is a tame extension. Every henselian field of residue characteristic 0
is a tame field.
\end{theorem}

An extension $(L|K,v)$ of valued fields is called \bfind{immediate} if
the canonical embeddings $vK\hookrightarrow vL$ and $Kv\hookrightarrow
Lv$ are onto. A valued field is called \bfind{algebraically maximal} if
it does not admit proper immediate algebraic extensions; it is called
\bfind{separable-algebraically maximal} if it does not admit proper
immediate separable-algebraic extensions.

Take a valued field $(K,v)$ and denote the characteristic exponent of
$Kv$ by $p$. Then $(K,v)$ is a \bfind{Kaplansky field} if $vK$ is
$p$-divisible and $Kv$ does not admit any finite extension whose degree
is divisible by $p$.
%
%
All algebraically maximal Kaplansky fields are tame fields (cf.\
Corollary~\ref{amKt} below). But the converse does not hold since for a
tame field it is admissible that its residue field has finite separable
extensions with degree divisible by $p$. It is because of this fact that
the uniqueness of maximal immediate extensions will in general fail
(cf.\ \cite{[K--P--R]}). This is what makes the proof of model theoretic
results for tame fields so much harder than for algebraically maximal
Kaplansky fields.

\parm
Let us now give a quick survey on the basic results in the model theory
of valued fields that will lead up to the questions we will ask for tame
fields.

We take ${\cal L}_{\rm VF} =\{+,-,\cdot\,,\mbox{ }^{-1},0,1, {\cal O}\}$
to be the language of valued fields, where ${\cal O}$ is a binary
relation symbol for valuation divisibility. That is, ${\cal O} (a,b)$
will be interpreted by $va\geq vb$, or equivalently, $a/b$ being an
element of the valuation ring ${\cal O}_K\,$. We will write ${\cal O}
(x)$ in place of ${\cal O} (x,1)$ (note that ${\cal O}(a,1)$ says that
$va\geq v1=0$, i.e., $a\in {\cal O}_K$).

For $(K,v)$ and $(L,v)$ to be elementarily equivalent in the language of
valued fields, it is necessary that $vK$ and $vL$ are elementarily
equivalent in the language ${\cal L}_{\rm OG}=$
\linebreak
$\{+,-,0,<\}$ of ordered
groups, and that $Kv$ and $Lv$ are elementarily equivalent in the
language ${\cal L}_{\rm F} =\{+,-,\cdot\,,\mbox{ }^{-1},0,1\}$ of fields
(or in the language ${\cal L}_{\rm R}=\{+,-,\cdot\,,0,1\}$ of rings).
This is because elementary sentences about the value group and about the
residue field can be encoded in the valued field itself.

A main goal in this paper is to find additional conditions on
$(K,v)$ and $(L,v)$ under which these necessary conditions are also
sufficient, i.e., the following \bfind{Ax--Kochen--Ershov Principle}
(in short: AKE$^\equiv$ Principle) holds:

\begin{equation}                            \label{AKEequiv}
vK\equiv vL\>\wedge\> Kv\equiv Lv \;\;\;\Longrightarrow\;\;\;
(K,v)\equiv \mbox{$(L,v)$}\;.
\end{equation}

An AKE$^\prec$ Principle is the following analogue for elementary
extensions:
\begin{equation}                            \label{AKEprec}
(K,v)\subseteq (L,v)\>\wedge\> vK\prec vL\>\wedge\> Kv\prec Lv
\;\;\;\Longrightarrow\;\;\;
(K,v)\prec \mbox{$(L,v)$}\;.
\end{equation}

If $\cal M$ is an ${\cal L}$-structure and ${\cal M}'$ a substructure of
$\cal M$, then we will say that ${\cal M}'$ is existentially closed in
${\cal M}$ and write ${\cal M}'\ec {\cal M}$ if every existential
${\cal L}$-sentence with parameters from ${\cal M}'$ which holds in
${\cal M}$ also holds in ${\cal M}'$. For the meaning of ``existentially
closed in'' in the setting of fields and of ordered abelian
groups, see \cite{[K--P]}. Inspired by Robinson's Test, our basic
approach will
be to ask for criteria for a valued field to be existentially closed in
a given extension field. Replacing $\prec$ by $\ec$, we thus look for
conditions which ensure that the following AKE$^\exists$ Principle
holds:
\begin{equation}                            \label{AKE}
(K,v)\subseteq (L,v)\>\wedge\> vK\ec vL\>\wedge\> Kv\ec Lv
\;\;\;\Longrightarrow\;\;\; (K,v)\ec\mbox{$(L,v)$}\;.
\end{equation}
The conditions
\begin{equation}                   \label{sico}
vK\ec vL\mbox{\ \ and\ \ } Kv\ec Lv
\end{equation}
will be called \bfind{side conditions}. It is an easy exercise in model
theoretic algebra to show that these conditions imply that $vL/vK$ is
torsion free and that $Lv|Kv$ is \bfind{regular}, i.e., the
algebraic closure of $Kv$ is linearly disjoint from $Lv$ over $Kv$, or
equivalently, $Kv$ is relatively algebraically closed in $Lv$ and
$Lv|Kv$ is separable; cf.\ Lemma~\ref{sica}.

A valued field for which (\ref{AKE}) holds will be called an
\bfind{AKE$^\exists$-field}. A class {\bf C} of valued fields will be
called \bfind{AKE$^\equiv$-class} (or \bfind{AKE$^\prec$-class}) if
(\ref{AKEequiv}) (or (\ref{AKEprec}), respectively) holds for all
$(K,v),(L,v)\in$ {\bf C}, and it will be called
\bfind{AKE$^\exists$-class} if (\ref{AKE}) holds for all
$(K,v)\in$ {\bf C}. We will also say that {\bf C} is \bfind{relatively
complete} if it is an AKE$^\equiv$-class, and that {\bf C} is
\bfind{relatively model complete} if it is an AKE$^\prec$-class. Here,
``relatively'' means ``relative to the value groups and residue
fields''.

%
%

The following elementary classes of valued fields are known
to satisfy all or some of the above AKE Principles:

\sn
a) Algebraically closed valued fields satisfy all three AKE Principles.
They even admit quantifier elimination; this has been shown by Abraham
Robinson, cf.\ \cite{[Ro]}.

\sn
b) Henselian fields of residue characteristic $0$ satisfy all
three AKE Principles. These facts have been (explicitly or implicitly)
shown by Ax and Kochen \cite{[AK]} and independently by Ershov
\cite{[Er3]}. They admit quantifier elimination relative to their value
group and residue field, as shown in \cite{[D]}.

\sn
c) $p$-adically closed fields: again, these fields were treated by Ax
and Kochen \cite{[AK]} and independently by Ershov \cite{[Er3]}.

\sn
d) $\wp$-adically closed fields (i.e., finite extensions of $p$-adically
closed fields): for definitions and results see the monograph by
Prestel and Roquette \cite{[P--R]}.

\sn
e) Finitely ramified fields: this case is a generalization of c) and d).
These fields were treated by Ziegler \cite{[Zi]} and independently by
Ershov \cite{[Er5]}.

\sn
f) Algebraically maximal Kaplansky fields: again, these fields were
treated by Ziegler \cite{[Zi]} and independently by Ershov \cite{[Er4]}.

Every valued field admits a maximal immediate algebraic
extension and a maximal immediate extension. All of the above mentioned
valued fields have the common property that these extensions are
unique up to valuation preserving isomorphism. This has always been a
nice tool in the proofs of the model theoretic results. However, as we
will show in this paper, this uniqueness is not indispensible. In its
absence, one just has to work much harder.

We will show that tame fields form an
AKE$^\exists$-class, and we will prove further model theoretic results
for tame fields and separably tame fields.

\pars
In many applications (such as the proof of a Nullstellensatz), only
existential sentences play a role. In these cases, it suffices to have
an AKE$^\exists$ Principle at hand. There are situations where we cannot
even expect more than this principle. In order to present one, we will
need some definitions that will be fundamental for this paper.

Every finite extension $(L|K,v)$ of valued fields satisfies the
{\bf fundamental inequality}:
\begin{equation}                             \label{fundineq}
n\>\geq\>\sum_{i=1}^{\rm g} {\rm e}_i {\rm f}_i
\end{equation}
where $n=[L:K]$ is the degree
of the extension, $v_1,\ldots,v_{\rm g}$ are the distinct extensions of
$v$ from $K$ to $L$, ${\rm e}_i=(v_iL:vK)$ are the respective
ramification indices and ${\rm f}_i=[Lv_i:Kv]$ are the respective
inertia degrees. The extension is called \bfind{defectless} if
equality holds in (\ref{fundineq}). Note that ${\rm g}=1$ if $(K,v)$
is henselian, so the definition given in axiom (TE3) is a special case
of this definition.

A valued field $(K,v)$ is called \bfind{defectless} (or \bfind{stable})
if each of its finite extensions is defectless, and \bfind{separably
defectless} if each of its finite separable extensions is defectless. If
$\chara Kv=0$, then $(K,v)$ is defectless (this is a consequence of the
``Lemma of Ostrowski'', cf.\ (\ref{LoO}) below).

Now let $(L|K,v)$ be any extension of valued fields. Assume that $L|K$
has finite trans\-cendence degree. Then (by Corollary~\ref{fingentb}
below):
\begin{equation}                            \label{wtdgeq}
\trdeg L|K \>\geq\> \trdeg Lv|Kv \,+\, \dim_\Q \Q\otimes vL/vK\;.
\end{equation}
%
%
%
We will say that $(L|K,v)$ is \bfind{without transcendence defect} if
equality holds in (\ref{wtdgeq}). If $L|K$ does not have finite
trans\-cendence degree, then we will say that $(L|K,v)$ is without
transcendence defect if every subextension of finite trans\-cendence
degree is. In Section~\ref{sectmwtd} we will prove:

\begin{theorem}                                       \label{wtA}
Every extension without transcendence defect of a henselian defectless
field satisfies the AKE$^\exists$ Principle.
\end{theorem}
\n
Note that it is not in general true that an extension without
transcendence defect of a henselian defectless field will satisfy the
AKE$^\prec$ Principle. There are such extensions that satisfy the side
conditions, for which the lower field is algebraically closed while the
upper field is not even henselian. A different and particularly
interesting example where both the lower and the upper field are
henselian and defectless is given in Theorem 3 of \cite{[K5]}.


A valued field $(K,v)$ has \bfind{equal characteristic} if $\chara
K=\chara Kv$. The following is the main theorem of this paper:

\begin{theorem}                             \label{tameAKE}
The class of all tame fields is an AKE$^\exists$-class and an
AKE$^\prec$-class.
The class of all tame fields of equal characteristic is an
AKE$^\equiv$-class.
%
\end{theorem}
This theorem, originally proved in \cite{[K1]}, has been applied in
\cite{[K6]} to study the structure of the Zariski space of all places of
an algebraic function field in positive characteristic.

\pars
As an immediate consequence of the foregoing theorem,
we get the following criterion for decidability:
\begin{theorem}                               \label{dec}
Let $(K,v)$ be a tame field of equal characteristic. Assume that the
theories $\mbox{\rm Th}(vK)$ of its value group (as an ordered group)
and $\mbox{\rm Th}(Kv)$ of its residue field (as a field) both
admit recursive elementary axiomatizations. Then also the theory of
$(K,v)$ as a valued field admits a recursive elementary axiomatization
and is decidable.
\end{theorem}
Indeed, the axiomatization of $\mbox{\rm Th}(K,v)$ can be taken to
consist of the axioms of tame fields of equal characteristic $\chara K$,
together with the translations of the axioms of $\mbox{\rm Th}(vK)$ and
$\mbox{\rm Th}(Kv)$ to the language of valued fields (cf.\
Lemma~\ref{elprvgrf}).

As an application, we will prove Theorem~\ref{sptame} in
Section~\ref{sectmttf} which includes the following decidability result:
\begin{theorem}
Take $q=p^n$ for some prime $p$ and some $n\in\N$, and an ordered
abelian group $\Gamma$. Assume that $\Gamma$ is divisible or
elementarily equivalent to the $p$-divisible hull of $\Z$. Then the
elementary theory of the power series field $\F_q((t^\Gamma))$ with
coefficients in $\F_q$ and exponents in $\Gamma$, endowed with its
canonical valuation $v_t\,$, is decidable.
\end{theorem}

\parm
Here are our results for separably defectless and separably tame
fields, which we will prove in Section~\ref{sectstsd}:

\begin{theorem}                             \label{septameAKE}
a) \ Take a separable extension $(L|K,v)$ without transcendence defect
of a henselian separably defectless field such that $vK$ is cofinal in
$vL$. Then the extension satisfies the AKE$^\exists$ Principle.
\n
b) \ Every separable extension $(L|K,v)$ of a separably tame field
satisfies the AKE$^\exists$ Principle.
\end{theorem}
We do not know whether the cofinality condition can be dropped when
$(K,v)$ is not a separably tame field.


\pars
We will deduce our model theoretic results from two main theorems which
we originally proved in \cite{[K1]}. The first theorem is a
generalization of the ``Grauert--Remmert Stability Theorem''. It deals
with function fields $F|K$, i.e., $F$ is a finitely generated field
extension of $K$ (for our purposes it is not necessary to ask that the
transcendence degree is $\geq 1$). For the following theorem, see
\cite{[K10]}:

\begin{theorem}                           \label{ai}
Let $(F|K,v)$ be a valued function field without transcendence defect.
If $(K,v)$ is a defectless field, then also $(F,v)$ is a defectless
field.
\end{theorem}

In \cite{[K--K1]} we used Theorem~\ref{ai} to prove \bfind{elimination
of ramification} for valued function fields. Let us describe this result
as it is not only important for the proofs of our model theoretic
results, but also in valuation theoretical approaches to resolution of
singularities.

The \bfind{henselization} of a valued field $(L,v)$ will be denoted by
$(L,v)^h$ or simply $L^h$. It is ``the minimal'' extension of $(L,v)$
which is henselian; for details, see Section~\ref{sectprel}. The
henselization is an immediate separable-algebraic extension. Hence
every separable-algebraically maximal valued field is henselian.

We will say that a valued function field $(F|K,v)$ is \bfind{strongly
inertially generated}, if there is a transcendence basis
\[
{\cal T}= \{x_1, \ldots,x_r,y_1,\ldots, y_s\}
\]
of $(F|K,v)$ such that
\pars
a) \ $vF= vK({\cal T})=vK\oplus\Z vx_1\oplus \ldots\oplus
\Z vx_r$,\par
b) \ $y_1v,\ldots,y_sv$ form a separating transcendence
basis of $Fv|Kv$,
\sn
and there is some $a$ in some henselization $F^h$ of $(F,v)$ such that
$F^h=K({\cal T})^h(a)$, $va=0$ and $K({\cal T})v(av)|K({\cal T})v$ is
separable of degree equal to $[K({\cal T})^h(a):K({\cal T})^h]$. The
latter means that $F$ lies in the ``absolute inertia field''
$K({\cal T})^i$ of $(K({\cal T}),v)$ (see definition in
Section~\ref{secttpw}).

The following is Theorem~3.4 of \cite{[K--K1]}:

\begin{theorem}                             \label{hrwtd}
Take a defectless field $(K,v)$ and a valued function field $(F|K,v)$
without transcendence defect. Assume that $Fv|Kv$ is a separable
extension and $vF/vK$ is torsion free. Then $(F|K,v)$ is strongly
inertially generated. In fact, for each transcendence basis ${\cal T}$
that satisfies conditions a) and b) there is an element $a$ with the
required properties in every henselization of $F$.
\end{theorem}

The second fundamental theorem, originally proved in \cite{[K1]}, is a
structure theorem for immediate function fields over tame or separably
tame fields (cf.\ \cite{[K2]}, \cite{[K11]}).

\begin{theorem}                \label{stt3}
Take an immediate function field $(F|K,v)$ of trans\-cendence degree 1.
Assume that $(K,v)$ is a tame field, or that $(K,v)$ is a separably
tame field and $F|K$ is separable. Then
\begin{equation}
\mbox{there is $x\in F$ such that }\; (F^h,v)\,=\,(K(x)^h,v)\;.
\end{equation}
\end{theorem}
For valued fields of residue characteristic 0, the assertion is a direct
consequence of the fact that every such field is defectless (in fact,
every $x\in F\setminus K$ will then do the job). In contrast to this,
the case of positive residue characteristic requires a much deeper
structure theory of immediate extensions of valued fields,
in order to find suitable elements $x$.

Theorem~\ref{stt3} is also used in \cite{[K--K2]}. For a survey on a
valuation theoretical approach to resolution of singularities and its
relation to the model theory of valued fields, see \cite{[K4]}.

%
%
\section{Valuation theoretical preliminaries}         \label{sectprel}
%

%
%
\subsection{Some general facts}
We will denote the algebraic closure of a field $K$ by $\tilde{K}$.
Whenever we have a valuation $v$ on $K$, we will automatically fix an
extension of $v$ to the algebraic closure $\tilde{K}$ of $K$. It does
not play a role which extension we choose, except if $v$ is also given
on an extension field $L$ of $K$; in this case, we choose the extension
to $\tilde{K}$ to be the restriction of the extension to $\tilde{L}$. We
say that $v$ is \bfind{trivial} on $K$ if $vK=\{0\}$. If the valuation
$v$ of $L$ is trivial on the subfield $K$, then we may assume that $K$
is a subfield of $Lv$ and the residue map $K\ni a\mapsto av$ is the
identity.

We will denote by $K\sep$ the separable-algebraic closure of $K$, and by
$K^{1/p^{\infty}}$ its perfect hull. If $\Gamma$ is an ordered abelian
group and $p$ a prime, then we write $\frac{1}{p^{\infty}}\Gamma$ for
the $p$-divisible hull of $\Gamma$, endowed with the unique extension of
the ordering from $\Gamma$. We leave the easy proof of the following
lemma to the reader.

\begin{lemma}                               \label{KsacP}
If $K$ is an arbitrary field and $v$ is a valuation on $K\sep$, then
$vK\sep$ is the divisible hull of $vK$, and $(Kv)\sep
\subseteq K\sep v$. If in addition $v$ is nontrivial on $K$, then
$K\sep v$ is the algebraic closure of $Kv$.

Every valuation $v$ on $K$ has a unique extension to $K^{1/p^{\infty}}$,
and it satisfies $vK^{1/p^{\infty}}=\frac{1}{p^{\infty}}vK$ and
$K^{1/p^{\infty}}v= (Kv)^{1/p^{\infty}}$.
\end{lemma}

For the easy proof of the following lemma, see \cite{[B]}, chapter VI,
\S10.3, Theorem~1.
\begin{lemma}                                      \label{prelBour}
Let $(L|K,v)$ be an extension of valued fields. Take elements $x_i,y_j
\in L$, $i\in I$, $j\in J$, such that the values $vx_i\,$, $i\in I$,
are rationally independent over $vK$, and the residues $y_jv$, $j\in
J$, are algebraically independent over $Kv$. Then the elements
$x_i,y_j$, $i\in I$, $j\in J$, are algebraically independent over $K$.

Moreover, write
\begin{equation}                            \label{polBour}
f\>=\> \displaystyle\sum_{k}^{} c_{k}\,
\prod_{i\in I}^{} x_i^{\mu_{k,i}} \prod_{j\in J}^{} y_j^{\nu_{k,j}}\in
K[x_i,y_j\mid i\in I,j\in J]
\end{equation}
in such a way that whenever $k\ne\ell$, then
there is some $i$ s.t.\ $\mu_{k,i}\ne\mu_{\ell,i}$ or some $j$ s.t.\
$\nu_{k,j}\ne\nu_{\ell,j}\,$. Then
\begin{equation}                            \label{value}
vf\>=\>\min_k\, v\,c_k \prod_{i\in I}^{}
x_i^{\mu_{k,i}}\prod_{j\in J}^{} y_j^{\nu_{k,j}}\>=\>
\min_k\, vc_k\,+\,\sum_{i\in I}^{} \mu_{k,i} v x_i\;.
\end{equation}
That is, the value of the polynomial $f$ is equal to the least of the
values of its monomials. In particular, this implies:
\begin{eqnarray*}
vK(x_i,y_j\mid i\in I,j\in J) & = & vK\oplus\bigoplus_{i\in I}
\Z vx_i\\
K(x_i,y_j\mid i\in I,j\in J)v & = & Kv\,(y_jv\mid j\in J)\;.
\end{eqnarray*}
The valuation $v$ on $K(x_i,y_j\mid i\in I,j\in J)$ is uniquely
determined by its restriction to $K$, the values $vx_i$ and the fact
that the residues $y_jv$, $j\in J$, are algebraically independent over
$Kv$.

The residue map on $K(x_i,y_j\mid i\in I,j\in J)$ is uniquely determined
by its restriction to $K$, the residues $y_jv$, and the fact that values
$vx_i\,$, $i\in I$, are rationally independent over $vK$.
%
\end{lemma}

We give two applications of this lemma.
\begin{corollary}                              \label{fingentb}
Take a valued function field $(F|K,v)$ without transcendence defect and
set $r=\dim_\Q \Q\otimes vF/vK$ and $s=\trdeg Fv|Kv$. Choose elements
$x_1,\ldots,x_r$, $y_1,\ldots,y_s\in F$ such that the values $vx_1,
\ldots,vx_r$ are rationally independent over $vK$ and the residues
$y_1v,\ldots,y_sv$ are algebraically independent over $Kv$. Then
${\cal T}=\{x_1,\ldots,x_r,y_1,\ldots,y_s\}$ is a transcendence basis
of $F|K$. Moreover, $vF/vK$ and the extension $Fv|Kv$ are finitely
generated.
\end{corollary}
\begin{proof}
By the foregoing theorem, the elements $x_1,\ldots,x_r,y_1,\ldots,
y_s$ are algebraically independent over $K$. Since $\trdeg F|K=r+s$ by
assumption, these elements form a transcendence basis of $F|K$.

It follows that the extension $F|K({\cal T})$ is finite. By the
fundamental inequality (\ref{fundineq}), this yields that $vF/
vK({\cal T})$ and $Fv| K({\cal T})v$ are finite. Since already
$vK({\cal T})/vK$ and $K({\cal T})v|Kv$ are finitely generated
by the foregoing lemma, it follows that also $vF/vK$ and
$Fv|Kv$ are finitely generated.
\end{proof}

If $(L|K,v)$ is an extension
of valued fields, then a transcendence basis ${\cal T}$ of $L|K$ will
be called a \bfind{standard valuation transcendence basis} of $(L,v)$
over $(K,v)$ if ${\cal T}=\{x_i,y_j\mid i\in I, j\in J\}$ where the
values $vx_i$, $i\in I$, form a maximal set of values in $vL$ rationally
independent over $vK$, and the residues $y_jv$, $j\in J$, form a
transcendence basis of $Lv|Kv$. Note that if $(L|K,v)$ is of finite
transcendence degree and admits a standard valuation transcendence
basis, then it is an extension without transcendence defect. Note also
that the transcendence basis ${\cal T}$ given in Theorem~\ref{hrwtd}
is a standard valuation transcendence basis.

\begin{corollary}                           \label{svtb->wtd}
If a valued field extension admits a standard valuation transcendence
basis, then it is an extension without transcendence defect.
\end{corollary}
\begin{proof}
Let $(L|K,v)$ be an extension with standard valuation transcendence basis
${\cal T}$, and $F|K$ a subextension of $L|K$ of finite transcendence
degree. We have to show that equality holds in (\ref{wtdgeq}) for $F$ in
place of $L$. Since $F|K$ is finitely generated, there is a finite
subset ${\cal T}_0\subseteq {\cal T}$ such that all generators
of $F$ are algebraic over $K({\cal T}_0)$. Then ${\cal T}_0$ is a
standard valuation transcendence basis of $(F({\cal T}_0)|K,v)$, and it
follows from Lemma~\ref{prelBour} that equality holds in (\ref{wtdgeq})
for $F':= F({\cal T}_0)$ in place of $L$. But as $\trdeg F'|K=\trdeg
F'|F+ \trdeg F|K$, $\trdeg F'v|Kv=\trdeg F'v|Fv+ \trdeg Fv|Kv$ and
$\dim_\Q \Q\otimes vF'/vK =\dim_\Q \Q\otimes vF'/vF +\dim_\Q \Q\otimes
vF/vK$, it follows that
\begin{eqnarray*}
\trdeg F'|K & = & \trdeg F'|F\,+\,\trdeg F|K \\
 & \geq & \trdeg F'v|Fv\,+\,\dim_\Q \Q\otimes vF'/vF\>+\>
\trdeg Fv|Kv\,+\,\dim_\Q \Q\otimes vF/vK\\
 & = & \trdeg F'v|Kv\,+\,\dim_\Q \Q\otimes vF'/vK\\
 & = & \trdeg F'|K\;,
\end{eqnarray*}
hence equality must hold. Since the inequality (\ref{wtdgeq}) holds for
the two extensions $(F'|F,v)$ and $(F|K,v)$, we find that $\trdeg F|K
=\trdeg Fv|Kv+\dim_\Q \Q\otimes vF/vK$ must hold.
\end{proof}


Every valued field $(L,v)$ admits a \bfind{henselization}, that is, a
minimal algebraic extension which is henselian. All henselizations are
isomorphic over $L$, so we will frequently talk of {\it the}
henselization of $(L,v)$, denoted by $(L,v)^h$, or simply $L^h$. The
henselization becomes unique in absolute terms once we fix an extension
of the valuation $v$ from $L$ to its algebraic closure. All
henselizations are immediate separable-algebraic extensions. They are
minimal henselian extensions of $(L,v)$ in the following sense: if
$(F,v')$ is a henselian extension field of $(L,v)$, then there is a
unique embedding of $(L^h,v)$ in $(F,v')$. This is the \bfind{universal
property of the henselization}. We note that every algebraic extension
of a henselian field is again henselian.

%
%
\subsection{The defect}                     \label{sectdef}
Assume that $(L|K,v)$ is a finite extension and the extension of $v$
from $K$ to $L$ is unique (which is always the case when $(K,v)$ is
henselian). Then the Lemma of Ostrowski (cf.\ \cite{[En]}, \cite{[Ri]},
\cite{[K2]}) says that
\begin{equation}                            \label{LoO}
[L:K]\;=\; (vL:vK)\cdot [Lv:Kv]\cdot p^\nu \;\;\;\mbox{ with }\nu\geq 0
\end{equation}
where $p$ is the characteristic exponent of $Kv$. The factor
\[\mbox{\rm d}(L|K,v)\>:=\> p^\nu \>=\> \frac{[L:K]}{(vL:vK)[Lv:Kv]}\]
is called the \bfind{defect} of the extension $(L|K,v)$. If $\nu>0$,
then we speak of a \bfind{nontrivial} defect. If $[L:K]=p$ then
$(L|K,v)$ has nontrivial defect if and only if it is immediate.
If $\mbox{\rm d}(L|K,v)=1$, then $(L|K,v)$ is a defectless
extension. Note that $(L|K,v)$ is always defectless if $\chara Kv=0$.

The following lemma shows that the defect is multiplicative.
This is a consequence of the multiplicativity of the degree of field
extensions and of ramification index and inertia degree. We leave the
straightforward proof to the reader.
\begin{lemma}                                       \label{md}
Take a valued field $(K,v)$. If $L|K$ and $M|L$ are finite extensions
and the extension of $v$ from $K$ to $M$ is unique, then
\begin{equation}         \label{pf}
\mbox{\rm d}(M|K,v) = \mbox{\rm d}(M|L,v)\cdot\mbox{\rm d}(L|K,v)
\end{equation}
In particular, $(M|K,v)$ is defectless if and only if $(M|L,v)$ and
$(L|K,v)$ are defectless.
\end{lemma}

The next lemma follows from Lemma~2.5 of \cite{[K8]}:

\begin{lemma}                               \label{idld}
Take an arbitrary immediate extension $(F|K,v)$ and an algebraic
extension $(L|K,v)$ of which every finite subextension admits a unique
extension of the valuation and is defectless. Then $F|K$ and $L|K$ are
linearly disjoint.
\end{lemma}

A valued field $(K,v)$ is called \bfind{inseparably defectless} if
equality holds in (\ref{fundineq}) for every finite purely inseparable
extension $L|K$. From the previous lemma, we obtain:

\begin{corollary}                     \label{iers}
Every immediate extension of a defectless field is regular. Every
immediate extension of an inseparably defectless field is separable.
\end{corollary}

The following is an important theorem, as passing to henselizations will
frequently facilitate our work.

\begin{theorem}                                    \label{dl-hdl}
Take a valued field $(K,v)$ and fix an extension of $v$ to $\tilde{K}$.
Then $(K,v)$ is defectless if and only if its henselization $(K,v)^h$ in
$(\tilde{K},v)$ is defectless. The same holds for ``separably
defectless'' and ``inseparably defectless'' in place of
``defectless''.
\end{theorem}
\begin{proof}
For ``separably defectless'', our assertion follows directly from
\cite{[En]}, Theorem (18.2). The proof of that theorem can easily be
adapted to prove the assertion for ``inseparably defectless'' and
``defectless''. See \cite{[K2]} for more details.
\end{proof}

Since a henselian field has a unique extension of the valuation to every
algebraic extension field, we obtain:

\begin{corollary}                               \label{lostr0}
Every valued field $(K,v)$ with $\chara Kv=0$ is a defectless field.
\end{corollary}

\begin{corollary}                                \label{mdc1}
A valued field $(K,v)$ is defectless if and only if $\mbox{\rm d}
(L|K^h,v)=1$ for every finite extension $L|K^h$.
\end{corollary}

Using this corollary together with Lemma~\ref{md}, one shows:
\begin{corollary}                                \label{mdc2}
Every finite extension of a defectless field is again a
defectless field.
\end{corollary}

\subsection{Tame and purely wild extensions}    \label{secttpw}
The interaction of tame extensions with the defect is described in the
following result, which is Proposition 2.8 of \cite{[K8]}:

\begin{proposition}                                  \label{dlta}
Let $(K,v)$ be a henselian field and $(N|K,v)$ an arbitrary tame
extension. If $L|K$ is a finite extension, then
\[
\mbox{\rm d}(L|K,v)\>=\> \mbox{\rm d}(L.N | N,v)\;.
\]
%
%
Hence, $(K,v)$ is a defectless field if and only if $(N,v)$ is a
defectless field. The same holds for ``separably defectless'' and
``inseparably defectless'' in place of ``defectless''.
\end{proposition}

We will denote by $K^r$ the ramification field of the normal extension
$(K\sep|K,v)$, and by $K^i$ its inertia field. As both fields contain
the decomposition field of $(K\sep|K,v)$, which is the henselization of
$K$ inside of $(K\sep,v)$, they are henselian.

The next lemma follows from general ramification theory; see
\cite{[En]}, \cite{[K2]}.
\begin{lemma}                               \label{Krmte}
Take a henselian field $(K,v)$.
\sn
a) \ If $(L|K,v)$ is an algebraic extension and $L'$ an intermediate
field, then $(L|K,v)$ is tame if and only if $(L'|K,v)$ and $(L|L',v)$
are.
\sn
b) \ The field $K^r$ is the unique maximal tame extension of $(K,v)$,
and $(K^r)^r=K^r$.
\end{lemma}

An algebraic extension of a henselian field is called \bfind{purely
wild} if it is linearly disjoint from every tame extension.
We will call $(K,v)$ a \bfind{purely wild field} if
$(\tilde{K}|K,v)$ is a purely wild extension.

Lemma~\ref{idld} immediately yields important examples of purely wild
extensions:

\begin{corollary}                           \label{iepw}
Every immediate algebraic extension of a henselian field is purely wild.
\end{corollary}

Part b) of Lemma~\ref{Krmte} yields the following facts:
\begin{lemma}                               \label{Krmte2}
Take a henselian field $(K,v)$. An algebraic extension of $(K,v)$ is
purely wild if and only if it is linearly disjoint from $K^r$. Further,
$K^r$ is a purely wild field.
\end{lemma}
\sn

Since $K^r|K$ is by definition a separable extension, Lemmas~\ref{Krmte}
and~\ref{Krmte2} yield:

\begin{corollary}                           \label{ts}
Every tame extension of a henselian field is separable. Every purely
inseparable algebraic extension of a henselian field is purely wild.
\end{corollary}

From Lemma~\ref{Krmte}, one easily deduces part a) of the next lemma.
Part b) follows from the fact that $L^r=L.K^r$ for every algebraic
extension $L|K$.

\begin{lemma}                               \label{chartafi}
a)\ \ Let $(K,v)$ be a henselian field. Then $(K,v)$ is a tame field
if and only if $K^r=\tilde{K}$. Similarly, $(K,v)$ is a separably tame
field if and only if $K^r=K\sep$. Further, $(K,v)$ is a purely wild
field if and only if $K^r=K$.
\sn
b)\ \ Every algebraic extension of a tame (or separably tame, or
purely wild, respectively) field is again a tame (or separably tame,
or purely wild, respectively) field.
\end{lemma}
%

The following theorem was proved by M.~Pank; cf.\ \cite{[K--P--R]}.

\begin{theorem}                                     \label{P}
Let $(K,v)$ be a henselian field with residue characteristic $p>0$.
There exist field complements $W_s$ of $K^r$ in $K\sep$ over $K$, i.e.,
$K^r.W_s = K\sep$ and $W_s$ is linearly disjoint from $K^r$ over $K$.
The perfect hull $W=W_s^{1/p^{\infty}}$ is a field complement of $K^r$
over $K$, i.e., $K^r.W = \tilde{K}$ and $W$ is linearly disjoint from
$K^r$ over $K$. The valued fields $(W_s,v)$ can be characterized as the
maximal separable purely wild extensions of $(K,v)$, and the
valued fields $(W,v)$ are the maximal purely wild extensions
of $(K,v)$.

Moreover, $vW=vW_s$ is the $p$-divisible hull of $vK$, and $Wv$ is the
perfect hull of $Kv$; if $v$ is nontrivial, then $Wv=W_s v$.
%
\end{theorem}
\n
In \cite{[K--P--R]}, a condition for the uniqueness of these complements is
given and its relation to Kaplansky's hypothesis A and the uniqueness
of maximal immediate extensions is explained.

\pars
We will need the following characterization of purely wild extensions:

\begin{lemma}                               \label{charpw}
An algebraic extension $(L|K,v)$ of henselian fields of residue
characteristic $p>0$ is purely wild if and only if $vL/vK$ is a
$p$-group and $Lv|Kv$ is purely inseparable.
\end{lemma}
\begin{proof}
By Zorn's Lemma, every purely wild extension is contained in a maximal
one. So our assertions on $vL/vK$ and $Lv|Kv$ follow from the
corresponding assertions of Theorem~\ref{P} for $vW$ and $Wv$.

For the converse, assume that $(L|K,v)$ is an extension of henselian
fields of residue characteristic $p>0$ such that $vL/vK$ is a $p$-group
and $Lv|Kv$ is purely inseparable. We have to show that $L|K$ is
linearly disjoint from every tame extension $(F|K,v)$. Since every
tame extension is a union of finite tame extensions, it suffices to show
this under the assumption that $F|K$ is finite. Then $[F:K]=(vF:vK)
[Fv:Kv]$. Since $p$ does not divide $(vF:vK)$ and $vL/vK$ is a
$p$-group, it follows that $vF\cap vL=vK$. As $vF+vL\subseteq v(F.L)$,
we have that
\[
(v(F.L):vL)\>\geq\>((vF+vL):vL)\>=\>(vF:(vF\cap vL))\>=\>(vF:vK)\;.
\]
Since $Fv|Kv$ is separable and $Lv|Kv$ is purely inseparable, these
extensions are linearly disjoint. As $(Fv).(Lv)\subseteq (F.L)v$, we
have that
\[
[(F.L)v:Lv]\>\geq\>[(Fv).(Lv):Lv]\>=\>[Fv:Kv]\;.
\]
Now we compute:
\[
[F.L:L]\>\geq\>(v(F.L):vL)[(F.L)v:Lv]\>\geq\>(vF:vK)[Fv:Kv]\>=\>
[F:K]\>\geq\>[F.L:L]\;,
\]
hence equality holds everywhere. This shows that $L|K$ is
linearly disjoint from $F|K$.
\end{proof}

In conjunction with equation (\ref{LoO}), this lemma shows:
\begin{corollary}                           \label{dpw=pp}
The degree of a finite purely wild extension $(L|K,v)$ of henselian
fields of residue characteristic $p>0$ is a power of $p$.
\end{corollary}

From Lemma~\ref{charpw} one also easily derives:

\begin{corollary}                           \label{wildtrans}
Take a henselian field $(K,v)$, an algebraic
extension $(L|K,v)$ and an intermediate field $L'$. Then $(L|K,v)$
is a purely wild extension if and only if $(L'|K,v)$ and $(L|L',v)$ are.
\end{corollary}

\pars
We use Proposition~\ref{dlta} and Theorem~\ref{P} to give the following
characterizations of defectless fields:

\begin{theorem}
Take a henselian field $(K,v)$. Then the following statements are
equivalent.
\sn
1) \ $(K,v)$ is a defectless field.
\sn
2) \ For some (or every) tame extension $(N|K,v)$, $(N,v)$ is a
defectless field.
\sn
3) \ $(K^r,v)$ is a defectless field.
\sn
4) \ Every finite purely wild extension of $(K,v)$ is defectless.
\sn
5) \ Every maximal purely wild extension of $(K,v)$ is defectless.
\pars
The same holds if ``defectless field'' is replaced by ``separably
defectless field'' in 1), 2) and 3) and ``purely wild extension''
is replaced by ``separable purely wild extension'' in 4) and 5).
\end{theorem}
\begin{proof}
The equivalence of 1), 2) and 3) both for ``defectless'' and ``separably
defectless'' follows from Proposition~\ref{dlta} and part b) of
Lemma~\ref{Krmte}. Similarly, the equivalence
of 4) and 5) for both properties follows from their definition for
arbitrary algebraic extensions. It suffices now to show the implication
4)$\Rightarrow$1), as the converse is trivial.

By Theorem~\ref{P}, there exists a field complement $W_s$ of $K^r$ over
$K$ in $K\sep$, and $W_s^{1/p^{\infty}}$ is a field complement of $K^r$
over $K$ in $\tilde{K}$. Consequently, given any finite extension (or
finite separable extension, respectively) $(L|K,v)$, there is a finite
subextension $N|K$ of $K^r|K$ and a finite subextension (or finite
separable subextension, respectively) $W_0|K$ of $W_s^{1/p^{\infty}}|K$
(or of $W_s|K$, respectively) such that $L\subseteq N.W_0$. It follows
that $N.L\subseteq N.W_0$. Since $(N|K,v)$ is a tame extension,
Lemma~\ref{dlta} shows that $d(L|K,v)= d(N.L|N,v)$ and
$d(N.W_0|N,v)=d(W_0|K,v)$. So we can compute:
\[
d(L|K,v)\>=\>d(N.L|N,v)\>\leq\>d(N.W_0|N,v)\>=\>d(W_0|K,v)\;.
\]
Hence if $(W_0|K,v)$ is defectless, then so is $(L|K,v)$. This proves
the desired implication.
%
%
%
\end{proof}

\pars
To conclude this section, we will prove the following technical result
that we shall use later.

\begin{lemma}                               \label{ifcont}
Take a valued field $(F,v)$ and suppose that $E$ is a subfield of
$F$ on which $v$ is trivial. Then $E\sep\subset F^i$. Further, if
$Fv|Ev$ is algebraic, then $(F.E\sep)v=(Fv)\sep$.
\end{lemma}
\begin{proof}
Our assumption implies that the residue map induces an embedding of $E$
in $Fv$. By ramification theory (\cite{[En]}, \cite{[K2]}), $F^iv=(Fv)
\sep$. Thus, $(Ev)\sep \subseteq F^iv$. Using Hensel's Lemma, one shows
that the inverse of the isomorphism $E \ni a \mapsto av\in Ev$ can be
extended from $Ev$ to an embedding of $(Ev)\sep$ in $F^i$. Its image is
separable-algebraically closed and contains $E$. Hence, $E\sep\subset
F^i$. Further, $(F.E\sep)v$ contains $E\sep v$, which by
Lemma~\ref{KsacP} contains $(Ev)\sep$. As $F.E\sep|F$ is algebraic, so
is $(F.E\sep)v|Fv$. Therefore, if $Fv|Ev$ is algebraic, then
$(F.E\sep)v$ is algebraic over $(Ev)\sep$ and hence
separable-algebraically closed. Since $(F.E\sep)v\subseteq
F^iv=(Fv)\sep$, it follows that $(F.E\sep)v=(Fv)\sep$.
\end{proof}

%
%
\subsection{Algebraically maximal and separable-algebraically
maximal fields}
All algebraically maximal and all separable-algebraically maximal fields
are henselian because the henselization is an immediate
separable-algebraic extension and therefore these fields must coincide
with their henselization. Every henselian defectless field is
algebraically maximal. However, the converse is not true in general:
algebraically maximal fields need not be defectless (see Example~3.25 in
\cite{[K9]}). But we will see in Corollary~\ref{cortame} below that it
holds for perfect fields of positive characteristic. More
generally, in \cite{[K8]} it is shown that a valued field of positive
characteristic is henselian and defectless if and only if it is
algebraically maximal and inseparably defectless. Note that for a valued
field of residue characteristic 0, ``henselian'', ``algebraically
maximal'' and ``henselian defectless'' are equivalent.

We will assume the reader to be familiar with the theory of pseudo
Cauchy sequences as developed in \cite{[Ka]}. Recall that a pseudo
Cauchy sequence $(a_\nu)_{\nu<\lambda}$ in $(K,v)$ is \bfind{of
transcendental type} if it fixes the value of every polynomial $f\in
K[X]$, that is, $vf(a_\nu)$ is constant for all large enough
$\nu<\lambda$. See \cite{[Ka]} for the proof of the following theorem.

\begin{theorem}                                      \label{fix}
A valued field $(K,v)$ is algebraically maximal if and only if every
pseudo Cauchy sequence in $(K,v)$ without a limit in $K$ is of
transcendental type.
\end{theorem}

We will need the following characterizations of algebraically maximal
and separable-algebraically maximal fields; cf.\ Theorems 1.4, 1.6 and
1.8 of \cite{[K8]}.
\begin{theorem}                             \label{exmaxvfth}
The valued field $(K,v)$ is algebraically maximal if and only if it is
henselian and for every polynomial $f\in K[X]$,
\begin{equation}                            \label{exmaxvf}
\exists x\in K\, \forall y\in K:\> vf(x)\geq vf(y)\;.
\end{equation}
Similarly, $(K,v)$ is separable-algebraically maximal if and only if
(\ref{exmaxvf}) holds for every separable polynomial $f\in K[X]$.
\end{theorem}
%

%
%
\section{The algebra of tame and separably tame fields}
%

%
%
\subsection{Tame fields}       \label{sectatf}  
From the definition of tame fields and the fact that every tame
extension is separable (Corollary~\ref{ts}), we obtain:
\begin{lemma}                               \label{thdp}
Every tame field is henselian, defectless and perfect.
\end{lemma}

In general, infinite algebraic extensions of defectless fields need not
again be defectless fields. For example, $\F_p(t)^h$ is a defectless
field by Theorem~\ref{ai} in conjunction with Theorem~\ref{dl-hdl}, but
the perfect hull of $\F_p(t)^h$ is a henselian field admitting an
immediate extension generated by a root of the polynomial
$X^p-X-\frac{1}{t}$ (cf.\ Example~3.12 of \cite{[K9]}). However, from
Lemmas~\ref{chartafi} and~\ref{thdp} we can deduce that every algebraic
extension of a tame field is a defectless field.


\pars
We give some characterizations for tame fields:

\begin{theorem}                    \label{tame}
Take a henselian field $(K,v)$ and denote by $p$ the characteristic
exponent of $Kv$. The following assertions are equivalent:\n
1) \  $(K,v)$ is a tame field,\n
2) \  $K^r$ is algebraically closed,\n
3) \  every purely wild extension $(L|K,v)$ is trivial,\n
4) \  $(K,v)$ is algebraically maximal and closed under purely
wild extensions by $p$-th roots,\n
5)\ \  $(K,v)$ is algebraically maximal, $vK$ is $p$-divisible
and $Kv$ is perfect.
\end{theorem}
\begin{proof}
The equivalence of 1) and 2) was stated already in part a) of
Lemma~\ref{chartafi}.
\sn
2)$\Rightarrow$3): \ By definition, a purely wild extension of $(K,v)$
must be linearly disjoint from $K^r=\tilde{K}$, hence trivial.
\sn
3)$\Rightarrow$4): \ Suppose that $(K,v)$ has no proper purely wild
extension. Then in particular, it has no proper purely wild extension by
$p$-th roots. From Corollary~\ref{iepw} we infer that $(K,v)$ admits no
proper immediate algebraic extensions, i.e., $(K,v)$ is algebraically
maximal.
\sn
4)$\Rightarrow$5): \
Assume now that $(K,v)$ is an algebraically maximal field closed under
purely wild extensions by $p$-th roots. Take $a\in K$. First, suppose
that $va$ is not divisible by $p$ in $vK$; then the extension $K(b)|K$
generated by an element $b\in\tilde{K}$ with $b^p = a$, together with
any extension of the valuation, satisfies $(vK(b):vK) \geq p =
[K(b):K]\geq (vK(b):vK)$. Hence, equality holds everywhere, and
(\ref{fundineq}) shows that $(vK(b):vK)=p$ and $K(b)v=Kv$. Hence by
Lemma~\ref{charpw}, $(K(b)|K,v)$ is purely wild, contrary to our
assumption on $(K,v)$. This shows that $vK$ is $p$-divisible.

Second, suppose that $va = 0$ and that $av$ has no $p$-th root in $Kv$.
Then $[K(b)v:Kv]\geq p=[K(b):K]\geq [K(b)v:Kv]$. Hence, equality holds
everywhere, and (\ref{fundineq}) shows that $vK(b)=vK$ and $[K(b)v:Kv]
=p$. It follows that $K(b)v|Kv$ is purely inseparable. Again by
Lemma~\ref{charpw}, the extension $(K(b)|K,v)$ is purely wild, contrary
to our assumption. This shows that $Kv$ is perfect.
\sn
5)$\Rightarrow$2): \
Suppose that $(K,v)$ is an algebraically maximal (and thus henselian)
field such that $vK$ is $p$-divisible and $Kv$ is perfect. Choose a
maximal purely wild extension $(W,v)$ in accordance to Theorem~\ref{P}.
Together with the last part of Theorem~\ref{P}, our condition on the
value group and the residue field yields that $(W|K,v)$ is immediate.
But since $(K,v)$ is assumed to be algebraically maximal, this extension
must be trivial. This shows that $\tilde{K}= K^r.W=K^r.K=K^r$.
\end{proof}

If the residue field $Kv$ does not admit any finite extension whose
degree is divisible by $p$, then in particular it must be perfect. Hence
we can deduce from the previous theorem:

\begin{corollary}                           \label{amKt}
Every algebraically maximal Kaplansky field is a tame field.
\end{corollary}

If $\chara Kv=0$, then $(K,v)$ is tame as soon as it is henselian, and
this is the case when it is algebraically maximal. If $\chara K=p>0$,
then every extension by $p$-th roots is purely inseparable and thus
purely wild. So the previous theorem together with Lemma~\ref{thdp}
yields:

\begin{corollary}                           \label{cortame}
a) \ A valued field $(K,v)$ of equal characteristic is tame if and only
if it is algebraically maximal and perfect.
\sn
b) \ If $(K,v)$ is an arbitrary valued field of equal characteristic,
then every maximal immediate algebraic extension of its perfect hull is
a tame field.
\sn
c) \ For perfect valued fields of equal characteristic, the
properties ``algebraically maximal'' and ``henselian and defectless''
are equivalent.
\end{corollary}

The implication 3)$\Rightarrow$1) of Theorem~\ref{tame} together
with Corollary~\ref{wildtrans} and Theorem~\ref{P} shows:
\begin{corollary}                           \label{Wtame}
Every complement $(W,v)$ as in Theorem~\ref{P} is a tame field.
\end{corollary}

\pars
The next corollary shows how to construct tame fields with suitable
prescribed value groups and residue fields.

\begin{corollary}                                      \label{Gk}
Take a perfect field $k$ of characteristic exponent $p$ and a
$p$-divisible ordered abelian group $\Gamma$. Then there exists a tame
field $K$ of characteristic exponent $p$ having $\Gamma$ as its value
group and $k$ as its residue field such that $K$ admits a standard
valuation transcendence basis over its prime field and the cardinality
of $K$ is equal to the maximum of the cardinalities of $\Gamma$ and $k$.
\end{corollary}
\begin{proof}
According to Theorem~2.14.\ of \cite{[K7]}, there is a valued field
$(K_0,v)$ of characteristic exponent $p$ with value group $\Gamma$ and
residue field $k$, and admitting a standard valuation transcendence
basis ${\cal T}$ over its prime field. Now take $(K,v)$ to be a maximal
immediate algebraic extension of $(K_0,v)$. Then $(K,v)$ is
algebraically maximal, and Theorem~\ref{tame} shows that it is a tame
field. Since it is an algebraic extension of $(K_0,v)$, it still admits
the same transcendence basis over its prime field.
%
%
If $v$ is trivial, then $\Gamma=\{0\}$ and $K=k$, whence $|K|=
\max\{|\Gamma|, |k|\}$. If $v$ is nontrivial, then $K$ and $\Gamma$ are
infinite and therefore, $|K|=\max\{\aleph_0,|{\cal T}|\}\leq
\max\{|\Gamma|,|k|\}\leq |K|$, whence again $|K|=\max\{|\Gamma|,|k|\}$.
\end{proof}

\parm
Now we will prove an important lemma on tame fields that we will need
in several instances.
\begin{lemma}                                     \label{trac}
Take a tame field $(L,v)$ and a relatively algebraically
closed subfield $K\subset L$. If in addition $Lv|Kv$ is an algebraic
extension, then $K$ is also a tame field and moreover, $vL/vK$ is
torsion free and $Kv = Lv$.
\end{lemma}
\begin{proof}
The following short and elegant version of the proof was given by
\ind{Florian Pop}.
Since $(L,v)$ is tame, it is henselian and perfect. Since $K$ is
relatively algebraically closed in $L$, it is henselian and perfect too.
Assume that $(K_1|K,v)$ is a finite purely wild extension; in view of
Theorem~\ref{tame}, we have to show that it is trivial. By
Corollary~\ref{dpw=pp}, the degree $[K_1:K]$ is a power of $p$, say
$p^m$. Since $K$ is perfect, $L|K$ and $K_1|K$ are separable extensions.
Since $K$ is relatively algebraically closed in $L$, we know that $L$
and $K_1$ are linearly disjoint over $K$. Thus, $K_1$ is relatively
algebraically closed in $K_1.L$, and
\[
[K_1.L:L] = [K_1:K] = p^m\;.
\]
Since $L$ is assumed to be a tame field, the extension $(K_1.L|L,v)$
must be tame. This implies that
\[(K_1.L)v\,|\,Lv\]
is a separable extension of degree $p^m$. By hypothesis, $Lv\,|\,Kv$ is
an algebraic extension, hence also $(K_1.L)v\,|\,Kv$ and
$(K_1.L)v\,|\,K_1v$ are algebraic. Furthermore, $(K_1.L,v)$ being a
henselian field and $K_1$ being relatively algebraically closed in
$K_1.L$, Hensel's Lemma shows that
\[
(K_1.L)v\,|\,K_1v
\]
must be purely inseparable. This yields that
\begin{eqnarray*}
p^m & = & [(K_1.L)v:Lv]_{\rm sep}^{} \leq
[(K_1.L)v:Kv]_{\rm sep}^{}
= [(K_1.L)v:K_1v]_{\rm sep}^{}\cdot
[K_1v:Kv]_{\rm sep}^{}\\
& = & [K_1v:Kv]_{\rm sep}^{}\leq [K_1v:Kv]
\leq [K_1:K] = p^m\;,
\end{eqnarray*}
showing that equality holds everywhere, which implies that
\[K_1v\,|\,Kv\]
is separable of degree $p^m$. Since $K_1|K$ was assumed to be purely
wild, we have $p^m = 1$ and the extension $K_1|K$ is trivial.

We have now shown that $K$ is a tame field; hence by Theorem~\ref{tame},
$vK$ is $p$-divisible and $Kv$ is perfect. Since $Lv|Kv$ is assumed to
be algebraic,
one can use Hensel's Lemma to show
that $Lv=Kv$ and that the torsion subgroup
of $vL/vK$ is a $p$-group. But as $vK$ is $p$-divisible, $vL/vK$ has
no $p$-torsion, showing that $vL/vK$ has no torsion at all.
\end{proof}

A similar fact holds for separably tame fields, as stated in
Lemma~\ref{Xsrac} below. Note that the conditions on the residue fields
is necessary, even if they are of characteristic 0 (cf.\ Example 3.9 in
\cite{[K7]}).

The following corollaries will show some nice properties of the class of
tame fields. They also admit generalizations to separably tame fields,
see Corollary~\ref{X+sff} below. First we show that a function field
over a tame field admits a so-called field of definition which is tame
and of finite rank, that is, its value group has only finitely many
convex subgroups. This is an important tool in the study of the
structure of such function fields.

\begin{corollary}                      \label{sff}
For every valued function field $F$ with given transcendence basis
${\cal T}$ over a tame field $K$, there exists a tame subfield $K_0$ of
$K$ of finite rank with $K_0v=Kv$ and $vK/vK_0$ torsion free, and
a function field $F_0$ with transcendence basis ${\cal T}$
over $K_0$ such that
\begin{equation}                       \label{sff1}
F= K.F_0
\end{equation}
and
\begin{equation}                       \label{sff2}
[F_0:K_0({\cal T})] = [F:K({\cal T})]\;.
\end{equation}
\end{corollary}
\begin{proof}
Let $F=K({\cal T})(a_1,\ldots,a_n)$. There exists a finitely generated
subfield $K_1$ of $K$ such that $a_1,\ldots,a_n$ are algebraic over
$K_1({\cal T})$ and $[F:K({\cal T})]=[K_1({\cal T})(a_1,\ldots,a_n):K_1
({\cal T})]$. This will still hold if we replace $K_1$ by any
extension field of $K_1$ within $K$. As a finitely generated field,
$(K_1,v)$ has finite rank. Now let
$y_j$, $j\in J$, be a system of elements in $K$ such that the residues
$y_j v$, $j\in J$, form a transcendence basis of $Kv$ over $K_1v$.
According to Lemma~\ref{prelBour}, the field $K_1(y_j|j\in J)$ has
residue field $K_1v(y_j v|j\in J)$ and the same value group as $K_1$,
hence it is again a field of finite rank. Let $K_0$ be the relative
algebraic closure of this field within $K$. Since by construction,
$Kv|K_1v(y_j v|j\in J)$ and thus also $Kv|K_0 v$ are algebraic, we can
infer from the preceding lemma that $K_0$ is a tame field with $K_0v=Kv$
and $vK/vK_0$ torsion free. As an algebraic extension of a field of finite
rank, it is itself of finite rank. Finally, the function field $F_0 =
K_0({\cal T})(a_1,\ldots,a_n)$ has transcendence basis ${\cal T}$ over
$K_0$ and satisfies equations (\ref{sff1}) and (\ref{sff2}).
\end{proof}

\begin{corollary}                \label{tint}
For every extension $(L|K,v)$ with $(L,v)$ a tame field, there exists a
tame intermediate field $L_0$ such that the extension $(L_0|K,v)$ admits
a standard valuation transcendence basis and the extension $(L|L_0,v)$
is immediate.
\end{corollary}
\begin{proof}
Take ${\cal T}=\{x_i,y_j\mid i\in I, j\in J\}$ where the values $vx_i$,
$i\in I$, form a maximal set of values in $vL$ rationally independent
over $vK$, and the residues $y_jv$, $j\in J$, form a transcendence basis
of $Lv|Kv$. With this choice, $vL/vK({\cal T})$ is a torsion group and
$Lv|K({\cal T})v$ is algebraic. Let $L_0$ be the relative
algebraic closure of $K({\cal T})$ within $L$. Then by Lemma~\ref{trac},
we have that $(L_0,v)$ is a tame field, that $Lv = L_0v$ and that
$vL/vL_0$ is torsion free and thus, $vL_0 = vL$. This shows that the
extension $(L|L_0,v)$ is immediate. On the other hand, ${\cal T}$ is a
standard valuation transcendence basis of $(L_0|K,v)$ by construction.
\end{proof}

%
%
\subsection{Separably tame fields}          \label{sectastf}
%
%
Note that separably tame fields of characteristic $0$ are tame and have
hence been covered in the previous section. So in this section we will
concentrate on valued fields of positive characteristic. Note also that
every trivially valued field is separably tame.

Since every finite separable-algebraic extension of a separably tame
field is a tame and thus defectless extension, a separably tame field is
always henselian and separably defectless. The converse is not true; it
needs additional assumptions on the value group and the residue field.
Under the assumptions that we are going to use frequently, the converse
will even hold for ``separable-algebraically maximal'' in place of
``henselian and separably defectless''. (Note that ``henselian and
separably defectless'' implies ``separable-algebraically maximal''.)
%

%
%

An \bfind{Artin-Schreier extension} of a field $K$ of characteristic
$p>0$ is an extension of degree $p$ generated by a root of a polynomial
of the form $X^p-X-a$ with $a\in K$. It is a Galois extension, and every
Galois extension of degree $p$ of a field of characteristic $p$ is an
Artin-Schreier extension.

\begin{theorem}    \label{septame}
Take a nontrivially valued field $(K,v)$ of characteristic $p>0$. The
following assertions are equivalent:
\n
1) \  $(K,v)$ is a separably tame field,\n
2) \  $K^r$ is separable-algebraically closed.\n
3) \  every separable purely wild extension $(L|K,v)$ is trivial,\n
4) \  $(K,v)$ is separable-algebraically maximal and closed under
purely wild Artin-Schreier extensions,\n
5) \  $(K,v)$ is separable-algebraically maximal, $vK$ is $p$-divisible
and $Kv$ is perfect.
\end{theorem}
\begin{proof}
The equivalence of 1) and 2) was stated already in part a) of
Lemma~\ref{chartafi}.
\sn
2)$\Rightarrow$3): \ By definition, a separable purely wild extension of
$(K,v)$ must be linearly disjoint from $K^r=K\sep$, hence trivial.
\sn
3)$\Rightarrow$4): \ Suppose that $(K,v)$ has no proper separable purely
wild extensions. Then in particular, $(K,v)$ admits no purely
wild Artin-Schreier extensions. Furthermore, $(K,v)$ admits no proper
separable-algebraic immediate extensions, as they would be purely wild.
Consequently, $(K,v)$ is separable-algebraically maximal.
\sn
4)$\Rightarrow$5): \ If $(K,v)$ is closed under purely wild
Artin-Schreier extensions and $v$ is nontrivial, then $vK$ is
$p$-divisible and $Kv$ is perfect (cf.\ Corollary 2.17 of \cite{[K7]}).
\sn
5)$\Rightarrow$2): \ Suppose that $(K,v)$ is a separable-algebraically
maximal field such that $vK$ is $p$-divisible and $Kv$ is perfect. Then
in particular,
%
%
$(K,v)$ is henselian. Choose a maximal separable purely wild extension
$(W_s,v)$ in accordance to Theorem~\ref{P}. Our condition on the value
group and the residue field yields that $(W_s|K,v)$ is immediate. But
since $(K,v)$ is assumed to be separable-algebraically maximal, this
extension must be trivial. This shows that $K\sep=K^r.W_s=K^r.K=K^r$.
\end{proof}

As in the case of tame fields, we derive the following results:
\begin{corollary} \n
a) \ Every separable-algebraically maximal Kaplansky field is a
separably tame field.
\sn
b) \ Every complement $(W_s,v)$ as in Theorem~\ref{P} is a separably
tame field.
\end{corollary}

\pars
Suppose that $(K,v)$ separably tame. Choose $(W_s,v)$ according to
Theorem~\ref{P}. Then by condition 3) of the theorem above, the
extension $(W_s|K,v)$ must be trivial. This yields that
$(K^{1/p^{\infty}} ,v)$ is the unique maximal purely wild extension of
$(K,v)$. Further, $(K,v)$ also satisfies condition 4) of the theorem.
From
%
%
Corollary 4.6 of \cite{[K8]} it follows that $(K,v)$ is dense in
$(K^{1/p^{\infty}},v)$, i.e., $K^{1/p^{\infty}}$ lies in the completion
of $(K,v)$. This proves:

\begin{corollary}                        \label{Xsept}
If $(K,v)$ is separably tame, then the perfect hull $K^{1/p^{\infty}}$
of $K$ is the unique maximal purely wild extension of $(K,v)$ and lies
in the completion of $(K,v)$. In particular, every immediate algebraic
extension of a separably tame field $(K,v)$ is purely inseparable and
included in the completion of $(K,v)$.
%
\end{corollary}
%
%

\begin{lemma}               \label{XKsK}
$(K,v)$ is a separably tame field if and only if $(K^{1/p^{\infty}},v)$
is a tame field. Consequently, if $(K^{1/p^{\infty}},v)$ is a tame
field, then $(K,v)$ is dense in $(K^{1/p^{\infty}},v)$.
\end{lemma}
\begin{proof}
Suppose that $(K,v)$ is a separably tame field. Then by the maximality
stated in the previous corollary, $(K^{1/p^{\infty}},v)$ admits no
proper purely wild algebraic extensions. Hence by Theorem~\ref{tame},
$(K^{1/p^{\infty}},v)$ is a tame field.

For the converse, suppose that $(K^{1/p^{\infty}},v)$ is a tame field.
Observe that the extension $(K^{1/p^{\infty}}|K,v)$ is purely wild and
contained in every maximal purely wild extension of $(K,v)$.
Consequently, if $(K^{1/p^{\infty}},v)$ admits no purely wild
extension at all, then $(K^{1/p^{\infty}},v)$ is the unique maximal
purely wild extension of $(K,v)$. Then in view of Theorem~\ref{P},
$K^{1/p^{\infty}}$ must be a field complement for $K^r$ over $K$ in
$\tilde{K}$. This yields that $K^r=K\sep$, hence by part b) of
Lemma~\ref{Krmte}, $(K\sep|K,v)$ is a tame extension,
%
%
showing that $(K,v)$ is a separably tame field. By
Corollary~\ref{Xsept}, it follows that $(K,v)$ is dense in
$(K^{1/p^{\infty}},v)$.
\end{proof}


The following lemma describes the interesting behaviour of separably
tame fields under composition of places.

\begin{lemma}                   \label{Xstco}
Take a separably tame field $(K,v)$ of characteristic $p>0$ and let $P$
be the place associated with $v$. Assume that $P=P_1P_2P_3$ where $P_1$
is a coarsening of $P$, $P_2$ is a place on $KP_1$ and $P_3$ is a place
on $KP_1P_2\,$. Assume further that $P_2$ is nontrivial (but $P_1$ and
$P_3$ may be trivial). Then $(KP_1,P_2)$ is a separably tame field. If
also $P_1$ is nontrivial, then $(KP_1,P_2)$ is a tame field.
\end{lemma}
\begin{proof}
By Theorem~\ref{tame}, $vK$ is $p$-divisible. The same is then true for
$v_{P_2} (KP_1)$. We wish to show that the residue field $KP_1P_2$ is
perfect. Indeed, assume that this were not the case. Then
%
%
there is an Artin-Schreier extension of $(K,P_1P_2)$ which adjoins a
$p$-th root to the residue field $KP_1P_2$ (cf.\ Lemma~2.13 of
\cite{[K7]}). Since this residue field extension is purely inseparable,
the induced extension of the residue field $Kv= KP_1P_2P_3$ can not be
separable of degree $p$. This shows that the Artin-Schreier extension is
a separable purely wild extension of $(K,v)$, contrary to our
assumption on $(K,v)$.

By Theorem~\ref{septame}, $(K,P)$ is separable-algebraically maximal.
%
%
This yields that the same is true for $(K,P_1P_2)$; indeed, if
$(L|K,P_1P_2)$ is an immediate extension, then $LP_1P_2=KP_1P_2$,
whence $LP_1P_2P_3=KP_1P_2P_3$, showing that also $(L|K,P)$ is
immediate.

If $P_1$ is trivial (hence w.l.o.g.\ equal to the identity map), then
$(KP_1,P_2)=(K,P_1P_2)$ is separable-algebraically maximal, and it
follows from Theorem~\ref{septame} that $(KP_1,P_2)$ is a separably tame
field.

Now assume that $P_1$ is nontrivial.
%
%
Suppose that there is a nontrivial immediate algebraic extension of
$(KP_1,P_2)$. Choose an element $b\notin KP_1$ in this extension, and
let $g$ be its minimal polynomial. Choose a monic polynomial $f\in K[X]$
such that $fP_1=g$, and a root $a$ of $f$. Then there is an extension of
$P_1$ to $K(a)$ such that $aP_1=b$. It follows from the fundamental
inequality that $K(a)P_1=KP_1(b)$ and that $(K(a),P_1)$ and $(K,P_1)$
have the same value group. But as $(KP_1(b)|KP_1,P_2)$ is immediate, it
now follows that also $(K(a)|K,P_1P_2P_3)$ is immediate. Note that we
can always choose $f$ to be separable as we may add a summand $cX$ with
$v_{P_1} c>0$, which does not change the image of $f$ under $P_1$. In
this way, we obtain a contradiction to the fact that $(K,P)$ is
separable-algebraically maximal. We have thus shown that $(KP_1,P_2)$ is
an algebraically maximal field, and it follows from Theorem~\ref{tame}
that $(KP_1,P_2)$ is a tame field.
\end{proof}

The following is an analogue of Lemma~\ref{trac}.
\begin{lemma}             \label{Xsrac}
Let $(L,v)$ be a separably tame field and $K\subset L$ a
relatively algebraically closed subfield of $L$. If the residue field
extension $Lv|Kv$ is algebraic, then $(K,v)$ is also a separably
tame field and moreover, $vL/vK$ is torsion free and $Kv = Lv$.
\end{lemma}
\begin{proof}
Since $K$ is relatively algebraically closed in $L$, it follows that
also $K^{1/p^{\infty}}$ is relatively algebraically closed in
$L^{1/p^{\infty}}$.
%
%
Since $(L,v)$ is a separably tame field, $(L^{1/p^{\infty}},v)$ is a
tame field by Lemma~\ref{XKsK}. From this lemma we also know that $Lv=
L^{1/p^{\infty}}v$ and $vL= vL^{1/p^{\infty}}$. Our assumption on
$Lv\,|\,Kv$ yields that the extension $L^{1/p^{\infty}}v \,|\,
K^{1/p^{\infty}}v$ is algebraic. From Lemma~\ref{trac} we can now infer
that $(K^{1/p^{\infty}},v)$ is a tame field with $K^{1/p^{\infty}}v=
L^{1/p^{\infty}}v=Lv$ and $vL^{1/p^{\infty}}/vK^{1/p^{\infty}}=
vL/vK^{1/p^{\infty}}$ torsion free. Again by Lemma~\ref{XKsK}, $(K,v)$
is thus a separably tame field with $Kv=K^{1/p^{\infty}}v=Lv$ and $vL/vK
= vL/vK^{1/p^{\infty}}$ torsion free.
\end{proof}

\begin{corollary}                           \label{X+sff}
Corollary~\ref{sff} also holds for separably tame fields in place of
tame fields. More precisely, if $F|K$ is a separable extension, then
$F_0$ and $K_0$ can be chosen such that $F_0|K_0$
%
%
is a separable extension. Moreover, if $vK$ is cofinal in
$vF$ then it can also be assumed that $vK_0$ is cofinal in $vF_0$.
\end{corollary}
\begin{proof}
Since the proof of Corollary~\ref{sff} only involves Lemma~\ref{trac},
it can be adapted by use of Lemma~\ref{Xsrac}. The first additional
assertion can be shown  using the fact that the finitely generated
separable extension $F|K$ is separably generated. The
second additional assertion is seen as follows. If $vF$ admits a biggest
proper convex subgroup, then let $K_0$ contain a nonzero element whose
value does not lie in this subgroup. If $vF$ and thus also $vK$
does not admit a biggest proper convex subgroup, then first choose $F_0$
and $K_0$ as in the (generalized) proof of Corollary~\ref{sff}; since
$F_0$ has finite rank, there exists some element in $K$ whose value does
not lie in the convex hull of $vF_0$ in $vF$, and adding this element to
$K_0$ and $F_0$ will make $vK_0$ cofinal in $vF_0$.
\end{proof}

With the same proof as for Corollary~\ref{tint}, but using
Lemma~\ref{Xsrac} in place of Lemma~\ref{trac}, one shows:
\begin{corollary}                           \label{+tint}
Corollary~\ref{tint} also holds for separably tame fields in
place of tame fields.
\end{corollary}

%
%
\section{Model theoretic preliminaries}
%
%
%
%
%
We will now discuss the axiomatization of valued fields and some of
their important properties. A valuation $v$ on a field $K$ can be given
in several ways. We can take the \bfind{valuation divisibility relation}
and formalize it as a binary predicate $R_v$ which in every valued field
is to be interpreted as
\[
R_v (x,y)\;\Longleftrightarrow\;vx\geq vy\;.
\]
But we can also take the valuation ring and formalize it as a unary
predicate ${\cal O}$ which in every valued field $(K,v)$ is to be
interpreted as
\[
{\cal O}(x)\;\Longleftrightarrow\;x\in {\cal O}\;.
\]
This predicate can be defined from the valuation divisibility relation
by
\[{\cal O}(x)\;\leftrightarrow\; R_v(x,1)\;.
\]
If we are working in the language of fields (what we usually do), then
the valuation divisibility relation can be defined from the predicate
${\cal O}$ by
\[
R_v(x,y)\;\leftrightarrow\;(y\not=0\,\wedge\,
{\cal O}(xy^{-1}))\>\vee\>x=0\;,
\]
whereas in general, it can not be defined using ${\cal O}$ and the
language of rings without the use of quantifiers, as in
\[
R_v(x,y)\;\leftrightarrow\;(\exists z\> yz=1\,\wedge \,
{\cal O}(xz))\>\vee\>x=0\;.
\]
This fact is only of importance for questions of quantifier elimination,
and only if one has decided to work in the language of rings. Note that
two fields are equivalent in the language of rings if and only if they
are equivalent in the language of fields. A similar assertion holds for
valued fields in the respective languages, and it also holds for the
notions ``elementary extension'' and ``existentially closed in'' in
place of ``equivalent''.
%

\pars
We prefer to write ``$vx\geq vy$'' in place of ``$R_v(x,y)$''. For
convenience, we define the following relations:
\begin{eqnarray*}
vx>vy & \leftrightarrow & vx\geq vy\>\wedge \>\neg(vy\geq vx)\\
vx=vy & \leftrightarrow & vx\geq vy\>\wedge \>vy\geq vx\;.
\end{eqnarray*}
The definitions for the reversed relations $vx\leq vy$ and $vx<vy$ are
obvious.

\pars
We will work in the language ${\cal L}_{\rm VF}$ of valued fields as
introduced in the introduction. The \bfind{theory of valued fields} is
the theory of fields (in the language ${\cal L}_{\rm F}$) together with
the axioms
\begin{axiom}
\ax{(V$0$)} $(\forall y\> vx\geq vy)\;\Leftrightarrow\; x=0$
\ax{(VT)} $v(x-y)\geq vx\>\vee\>v(x-y)\geq vy$
\end{axiom}
and the axioms which state that the value group is an ordered
abelian group:
\begin{axiom}
\ax{(VV$\not\!\mbox{\rm R}$)} $\neg (vx<vx)$
\ax{(VVT)} $vx<vy\,\wedge \,vy<vz\>\Rightarrow\>vx<vz$
\ax{(VVC)} $vx<vy\,\vee\,vx=vy\,\vee\,vx>vy$
\ax{(VVG)} $vx<vy\>\Rightarrow\>vxz<vyz$
\end{axiom}
(the group axioms for the value group follow from the group axioms for
the multiplicative group of the field).
%
%
%

The following facts are well-known; the easy proofs are left to the
reader.
\begin{lemma}                                    \label{elprvgrf}
Take a valued field $(K,v)$.\sn
a)\ \ For every sentence $\varphi$ in the language of ordered groups
there is a sentence $\varphi'$ in the language of valued fields such
that for every valued field $(K,v)$, $\varphi$ holds in $vK$ if and only
if $\varphi'$ holds in $(K,v)$.\sn
b)\ \ For every sentence $\varphi$ in the language of rings there is a
sentence $\varphi'$ in the language of valued fields such that for every
valued field $(K,v)$, $\varphi$ holds in $Kv$ if and only if
$\varphi'$ holds in $(K,v)$.\sn
\end{lemma}

As immediate consequences of this lemma, we obtain:
\begin{corollary}                                   \label{vgrfequiv}
If $(K,v)$ and $(L,v)$ are valued fields such that $(K,v)\equiv (L,v)$
in the language of valued fields, then $vK\equiv vL$ in the language of
ordered groups, and $Kv\equiv Lv$ in the language of rings (and
thus also in the language of fields). The same holds with $\prec$ or
$\ec$ in place of $\equiv\,$.
\end{corollary}

\begin{corollary}                                   \label{vgrfsat}
If $(K,v)$ is $\kappa$-saturated, then so are $vK$ (in the language of
ordered groups) and $Kv$ (in the language of fields).
\end{corollary}

The property of being henselian is axiomatized by the following axiom
scheme:
\begin{axiom}
\ax{(HENS)} $vy\geq 0\>\wedge\>\bigwedge_{1\leq i\leq n} vx_i\geq 0
\>\wedge\> v(y^n+x_1y^{n-1}+\ldots +x_{n-1}y +x_n)>0\\
\hspace*{4.6cm}\wedge\>v(ny^{n-1}+(n-1)x_1y^{n-2}+\ldots+x_{n-1})=0\\
\Rightarrow\> \exists z\; v(y-z)>0\,\wedge\,
z^n+x_1z^{n-1}+\ldots +x_{n-1}z +x_n=0\hfill(n\in\N)\,. \hspace{.3cm}$
\end{axiom}
Here we use one of the forms of Hensel's Lemma
%
%
to characterize henselian fields (see \cite{[K2]} for an extensive
collection). In view of Theorem~\ref{exmaxvfth}, also the property of
being algebraically maximal is easily axiomatized by axiom scheme (HENS)
together with the following axiom scheme:
\begin{axiom}
\ax{(MAXP)} $\exists y\forall z:\>v(y^n+x_1y^{n-1}+\ldots+x_{n-1}y+x_n)
\geq v(z^n+x_1z^{n-1}+\ldots+x_{n-1}z+x_n)\\
\mbox{ }\hfill(n\in\N)\,.\hspace{2cm}$
\end{axiom}
By the same theorem, the property of being separable-algebraically
maximal is axiomatized by axiom scheme (HENS) together with a version of
axiom scheme (MAXP) restricted to separable polynomials. This is
obtained by adding sentences that state that the coefficient of at least
one power $y^i$ for $i>0$ not divisible by the characteristic of the
field is nonzero.

The following was proved by Delon \cite{[D]} and Ershov \cite{[Er7]}.
For the case of valued fields of positive characteristic, we give an
alternative proof in \cite{[K8]}.
\begin{lemma}                               \label{hdlele}
The property of being henselian and defectless is elementary.
\end{lemma}

%
%
\section{The AKE$^\exists$ Principle}  \label{sectdAKEP}
%
%
%
\subsection{Necessary conditions for the AKE$^\exists$ Principle}

\parm
In this section we discuss tools for the proof of AKE$^\exists$
Principles and ask for those properties that a valued field must
have if it is an AKE$^\exists$-field.

We will need a model theoretic tool which we will apply to valued fields
as well as value groups and residue fields. We consider a countable
language ${\cal L}$ and ${\cal L}$-structures ${\eu B}$ and ${\eu A}^*$
with a common substructure ${\eu A}$. We will say that $\sigma$ is an
\bfind{embedding of ${\eu B}$ in ${\eu A}^*$ over ${\eu A}$} if it is an
embedding of ${\eu B}$ in ${\eu A}^*$ that leaves the universe $A$ of
${\eu A}$ elementwise fixed.

In what follows we will use Lemma~5.2.1.\ of \cite{[C--K]}, which states
that if ${\eu A}^*$ is $|B|^+$-saturated and every existential sentence
that holds in ${\eu B}$ also holds in ${\eu A}^*$, then ${\eu B}$ embeds
in ${\eu A}^*$.

\begin{proposition}                               \label{ec}
Let ${\eu A}\subseteq {\eu B}$ and ${\eu A}\subseteq {\eu A}^*$ be
extensions of ${\cal L}$-structures. If\/ ${\eu B}$ embeds over
${\eu A}$ in ${\eu A}^*$ and if\/ ${\eu A}\ec {\eu A}^*$, then ${\eu A}
\ec {\eu B}$. Conversely, if ${\eu A}\ec {\eu B}$ holds and if
${\eu A}^*$ is $|B|^+$-saturated, then ${\eu B}$ embeds over
${\eu A}$ in ${\eu A}^*$.
\end{proposition}
\begin{proof}
Since ${\eu A}$ is a substructure of ${\eu B}$ and of ${\eu A}^*$, both
$({\eu B},A)$ and $({\eu A}^*,A)$ are ${\cal L}(A)$-structures.

Suppose that $\sigma$ is an embedding of ${\eu B}$ over ${\eu A}$ in
${\eu A}^*$. Then every ${\cal L}(A)$-sentence will hold in $({\eu B},
A)$ if and only if it holds in $(\sigma {\eu B},A)$ (because isomorphic
structures are equivalent). Every existential ${\cal L}(A)$-sentence
$\varphi$ which holds in $({\eu B}, A)$ will then also hold in $({\eu
A}^*,A)$ since ${\eu A}^*$ is an extension of $\sigma {\eu B}$. If in
addition ${\eu A}\ec {\eu A}^*$, then $\varphi$ will also hold in
$({\eu A},A)$. This proves our first assertion.

Now suppose that ${\eu A}\ec {\eu B}$. Then every existential ${\cal L}
(A)$-sentence which holds in $({\eu B},A)$ also holds in $({\eu A},A)$
and, as $({\eu A}^*,A)$ is an extension of $({\eu A},A)$, also in $({\eu
A}^*, A)$. Now assume in addition that ${\eu A}^*$ is $|B|^+$-saturated.
Since $|A|\leq |B|< |B|^+$, also $({\eu A}^*,A)$ is $|B|^+$-saturated.
Hence by the lemma cited above, $({\eu B},A)$ embeds in $({\eu A}^*,A)$,
i.e., ${\eu B}$ embeds in ${\eu A}^*$ over ${\eu A}$.
\end{proof}

If we have an extension ${\eu A}\subseteq {\eu B}$ of
${\cal L}$-structures and want to show that ${\eu A}\ec {\eu B}$,
then by our proposition it suffices to show that ${\eu B}$ embeds
over ${\eu A}$ in some elementary extension ${\eu A}^*$ of ${\eu A}$.
This is the motivation for {\bf embedding lemmas},\index{embedding
lemma} which will play an important role later in our paper. When we
look for such embeddings, we can use the following very helpful
principles:
%

\begin{lemma}                                        \label{embfingen}
Let ${\eu A}\subseteq {\eu B}$ be an extension of ${\cal L}$-structures.
\sn
a) \ ${\eu A}$ is existentially closed in ${\eu B}$ if and only if it is
existentially closed in every substructure of ${\eu B}$ which is
finitely generated over ${\eu A}$.

\sn
b) \ Assume that ${\eu A}^*$ is a $|B|^+$-saturated extension of
${\eu A}$. If every substructure of ${\eu B}$ which is finitely
generated over ${\eu A}$ embeds over ${\eu A}$ in ${\eu A}^*$,
then also ${\eu B}$ embeds over ${\eu A}$ in ${\eu A}^*$.
\end{lemma}
\begin{proof}
a): \ If ${\eu A}$ is existentially closed in ${\eu B}$ then it is also
existentially closed in every substructure of ${\eu B}$ that contains
${\eu A}$ because an existential sentence that holds in this substructure
also holds in ${\eu B}$.

Every existential sentence only talks about finitely many elements,
hence it holds in $({\eu B},A)$ if and only if it holds in $({\eu B}_0,
A)$ where ${\eu B}_0$ is the substructure of ${\eu B}$ generated over
${\eu A}$ by these finitely many elements. Hence if ${\eu A}$ is
existentially closed in every such substructure, then it is
existentially closed in ${\eu B}$.

\sn
b) \ By what we have stated in the proof of part a) it follows that if
every substructure of ${\eu B}$ which is finitely generated over
${\eu A}$ embeds over ${\eu A}$ in ${\eu A}^*$, then every existential
sentence that holds in ${\eu B}$ will also hold in some image in ${\eu
A}^*$ of such a substructure, and hence it will hold in ${\eu A}^*$.
%
%
Using Lemma~5.2.1.\ of \cite{[C--K]}, we obtain that ${\eu B}$ embeds in
${\eu A}^*$ over ${\eu A}$.
\end{proof}

\pars
We will also need the following well known facts (which were proved,
e.g., in L.~van den Dries' thesis).

\begin{lemma}                               \label{sica}
a) \ Take an extension $G|H$ of torsion free abelian groups. If $H$ is
existentially closed in $G$ in the language ${\cal L}_{\rm G}=\{+,-,0\}$
of groups, then $G/H$ is torsion free.
\sn
b) \ Take a field extension $L|K$. If $K$ is existentially closed in
$L$ in the language ${\cal L}_{\rm F}$ of fields (or in the language
${\cal L}_{\rm R}$ of rings), then $L|K$ is regular.
\end{lemma}

An immediate consequence of the AKE$^\exists$ Principle (\ref{AKE}) is
the following observation:

\begin{lemma}                               \label{AKEex-am}
Every AKE$^\exists$-field is algebraically maximal.
\end{lemma}
\begin{proof}
Take a valued field $(K,v)$ which admits an immediate algebraic
extension $(L,v)$. Then by Lemma~\ref{sica} b), $K$ is not existentially
closed in $L$. Hence, $(K,v)$ is not existentially closed in $(L,v)$.
But $vK=vL$ and $Kv=Lv$, so that the conditions $vK\ec vL$ and
$Kv\ec Lv$ hold. This shows that $(K,v)$ is not an AKE$^\exists$-field.
\end{proof}
In particular, this lemma shows that every AKE$^\exists$-field must be
henselian.

\pars
A special case of the AKE$^\exists$ Principle is given if an extension
$(L|K,v)$ is immediate. Then, the side conditions are trivially
satisfied. We conclude that an AKE$^\exists$-field must in particular be
existentially closed in every immediate extension $(L,v)$. (We have used
this idea already in the proof of the foregoing lemma.) We can exploit
this fact by taking $(M,v)$ to be a maximal immediate extension of
$(K,v)$, to see which properties of $(M,v)$ are inherited by $(K,v)$ if
$(K,v)\ec (M,v)$. We know that $(M,v)$ has strong structural properties:
every pseudo Cauchy sequence has a limit (cf.\ \cite{[Ka]}), and it is
spherically complete (cf.\ \cite{[K2]}).
%
%


Since $(M,v)$ must coincide with its henselization which is an immediate
extension, it is henselian. By Theorem~31.21 of \cite{[W]},
%
%
$(M,v)$ is also a defectless field. Nevertheless, if $(K,v)$ is
henselian of residue characteristic $0$, then $(K,v)\prec (M,v)$,
which means that the elementary properties of $(M,v)$ are not stronger
than those of $(K,v)$. For other classes of valued fields, the situation
can be very different. Let us prove that every AKE$^\exists$-field is
henselian and defectless:

\begin{lemma}                               \label{Mod1Alg3}
Let $(K,v)$ be a valued field and assume that there is some maximal
immediate extension $(M,v)$ of $(K,v)$ which satisfies $(K,v)\ec
(M,v)$. Then $(K,v)$ is henselian and defectless. In particular, every
AKE$^\exists$-field is henselian and defectless.
\end{lemma}
\begin{proof}
Let $(E|K,v)$ be an arbitrary finite extension. Working in the language
of valued fields augmented by an additional predicate for a subfield, we
take $(E|K,v)^*$ to be a $|M|^+$-saturated elementary extension of
$(E|K,v)$.
%
%
Then $(E^*,v^*)$ and $(K^*,v^*)$ are $|M|^+$-saturated elementary
extensions of $(E,v)$ and $(K,v)$ respectively. Since by assumption
$(K,v)$ is existentially closed in $(M,v)$, Proposition~\ref{ec} shows
that we can embed $(M,v)$ over $(K,v)$ in $(K^*,v^*)$. We identify it
with its image in $(K^*,v^*)$. Since $(E^*|K^*,v^*)$ is an elementary
extension of $(E|K,v)$, the extensions $E|K$ and $K^*|K$ are linearly
disjoint. Therefore, $n := [E:K]=[E.M:M]$.

We will prove that the extension $(E.M,v^*)|(E,v)$ is immediate. Since
$E.M|M$ is algebraic and $vM = vK$, we know from the fundamental
inequality (\ref{fundineq}) that $v^*(E.M)/vK$ and hence
also $v^*(E.M)/vE$ is a torsion group. For the same reason, $Mv = Kv$
yields that $(E.M)v^*|Kv$ and hence also $(E.M)v^*|Ev$ is algebraic. On
the other hand, since $(E^*,v^*)$ is an elementary extension of
\mbox{$(E,v)$} we know by Lemma~\ref{sica} that $v^*E^*/vE$ is torsion
free and that $Ev$ is relatively algebraically closed in $E^*v$.
Combining these facts, we get that
\[
v^*(E.M) \>=\> vE \mbox{\ \ and\ \ } (E.M)v^* \>=\> Ev\;,
\]
showing that $(E.M,v^*)|(E,v)$ is immediate, as contended.

Since $(M,v)$ is maximal, it is a henselian and defectless field, as we
have mentioned above. Consequently,
\[
[E:K] \>=\> n \>=\> [E.M:M] \>=\> (v^*(E.M):vM)\cdot [(E.M)v^*:Mv]
\>=\> (vE:vK)\cdot [Ev:Kv]\;,
\]
which shows that $(E|K,v)$ is defectless and that the extension
of the valuation $v$ from $K$ to $E$ is unique. Since $(E,v)$ was an
arbitrary finite extension of $(K,v)$, this shows that $(K,v)$ is a
henselian and defectless field.
\end{proof}

%
%
%
\subsection{Extensions without transcendence defect}  \label{sectmwtd}
Our first goal in this section is to prove Theorem~\ref{wtA}. Take a
henselian and defectless field $(K,v)$ and an extension $(L|K,v)$
without transcendence defect. We choose $(K^*,v^*)$ to be an
$|L|^+$-saturated elementary extension of $(K,v)$. Since ``henselian''
is an elementary property, $(K^*,v^*)$ is henselian like $(K,v)$.
Further, it follows from Corollary~\ref{vgrfsat} that $K^*v^*$ is an
$|Lv|^+$-saturated elementary extension of $Kv$ and that $v^*K^*$ is a
$|vL|^+$-saturated elementary extension of $vK$. Assume that the side
conditions $vK\ec vL$ and $Kv\ec Lv$ hold. Then by Proposition~\ref{ec},
there exist embeddings
\[
\rho:\>vL \longrightarrow v^*K^*
\]
over $vK$ and
\[
\sigma:\;Lv \longrightarrow K^*v^*
\]
over $Kv$. Here, the embeddings of value groups and residue fields are
understood to be monomorphisms of ordered groups and of fields,
respectively.

We wish to prove that $(K,v)\ec (L,v)$. By Proposition~\ref{ec}, this
can be achieved by showing the existence of an embedding
\[
\iota:\> (L,v) \longrightarrow (K^*,v^*)
\]
over $K$, i.e., an embedding of $L$ in $K^*$ over $K$ preserving the
valuation, that is,
\[
\forall x\in L:\; x\in {\cal O}_L\Longleftrightarrow
\iota x \in {\cal O}_{K^*}\;.
\]

According to part b) of Lemma~\ref{embfingen}, such an embedding exists
already if it exists for every finitely generated subextension $(F|K,v)$
of $(L|K,v)$. In this way, we reduce our embedding problem to an
embedding problem for valued algebraic function fields $(F|K,v)$. Since
in the present case, $(L|K,v)$ is assumed to be an extension without
transcendence defect, the same holds for every finitely generated
subextension $(F|K,v)$. The case of such valued function fields is
covered by the following embedding lemma.

For a polynomial $f\in {\cal O}_K[X]$, we denote by $fv$ the polynomial
in $Kv[X]$ that is obtained from $f$ by replacing all its coefficients
by their residues.

\begin{lemma}                        \label{ael}
{\bf (Embedding Lemma I)}\n
Let $(F|K,v)$ a strongly inertially generated function field
and $(K^*,v^*)$ a henselian extension of $(K,v)$. Assume
that $vF/vK$ is torsion free and that $Fv|Kv$ is separable. If $\rho:\>
vF \longrightarrow v^*K^*$ is an embedding over $vK$ and $\sigma:\>Fv
\longrightarrow K^*v^*$ is an embedding over $Kv$, then there exists an
embedding $\iota:\>(F,v) \longrightarrow (K^*,v^*)$ over $(K,v)$ that
respects $\rho$ and $\sigma$, i.e., $v^*(\iota a) = \rho (va)$ and
$(\iota a)v^*= \sigma(av)$ for all $a \in F$.
\end{lemma}
\begin{proof}              
We choose a transcendence basis ${\cal T}$ and an element $a$ as in the
definition of an inertially generated function field.
First we will construct the embedding for $K({\cal T})$ and then we will
show how to extend it to $F$.

We choose elements $x_1',\ldots ,x_r'\in K^*$ such that
$v^* x_i' = \rho (vx_i)$, $1\leq i \leq r$. The values $v^*
x_1',\ldots ,v^* x_r'$ are rationally independent over $vK$
since the same holds for their preimages $vx_1,\ldots ,vx_r$ and
this property is preserved by every monomorphism over $vK$. Next, we
choose elements $y_1',\ldots ,y_s'\in {\cal O}_{K^*}^{\times}$ such
that $y_j'v^* = \sigma (y_j v)$, $1\leq j\leq s$. The residues
$y_1'v^* ,\ldots ,y_s'v^*$ are algebraically independent over
$Kv$ since the same holds for their preimages $y_1 v,\ldots,y_s v$ and
this property is preserved by every monomorphism over $Kv$.
Consequently, the elements $x_1',\ldots ,x_r'$ and $y_1',\ldots,y_s'$ as
well as the elements $x_1,\ldots ,x_r$ and $y_1,\ldots,y_s$ satisfy the
conditions of Lemma~\ref{prelBour}. Hence, both sets ${\cal T}$ and
${\cal T}'=\{x_1',\ldots ,x_r',y_1',\ldots,y_s'\}$ are algebraically
independent over $K$, so that the assignment
\[
x_i\mapsto x_i'\;,\;y_j\mapsto y_j'\;\;1\leq i\leq r\,,\; 1\leq j\leq s
\]
induces an isomorphism $\iota: K({\cal T}) \longrightarrow
K({\cal T}')$. Furthermore, for every $f\in K[{\cal T}]$, written
as in (\ref{polBour}),

\begin{eqnarray*}
v^*(\iota f) & =& \min_k\left(v^*c_k\,+\,\sum_{1\leq i\leq r} \mu_{k,i}
v^*x'_i\right) \>=\> \min_k\left(vc_k\,+\,\sum_{1\leq i\leq r} \mu_{k,i}
\rho vx_i\right) \\
 & =& \rho\, \min_k\left(vc_k\,+\,\sum_{1\leq i\leq r} \mu_{k,i}
vx_i\right)\>=\>\rho (vf)\;,
\end{eqnarray*}
showing that $\iota$  respects the restriction of $\rho$ to
$vK({\cal T})$. If $vf=0$, then
\[
fv\>=\>\left(\displaystyle\sum_{\ell}^{} c_{\ell}\,
\prod_{1\leq j\leq s} y_j^{\nu_{\ell,j}}\right)v
\>=\>\displaystyle\sum_{\ell}^{} (c_{\ell}v)
\prod_{1\leq j\leq s} (y_j v)^{\nu_{\ell,j}}
\]
where the sum runs only over those $\ell=k$ for which $\mu_{k,i}=0$ for
all $i$, and a similar formula holds for $(\iota f)v$ with the same
indices $\ell$. Hence,
\begin{eqnarray*}
(\iota f)v^* & =& \displaystyle\sum_{\ell}^{} (c_{\ell}v^*)
\prod_{1\leq j\leq s} (y_j v^*)^{\nu_{\ell,j}} \>=\>
 \displaystyle\sum_{\ell}^{} (c_{\ell}v)
\prod_{1\leq j\leq s} \sigma (y_j v)^{\nu_{\ell,j}} \\
 & =& \sigma\left(\displaystyle\sum_{\ell}^{} (c_{\ell}v)
\prod_{1\leq j\leq s}^{} (y_j v)^{\nu_{\ell,j}}\right) \>=\>\sigma
(fv)\;,
\end{eqnarray*}
showing that $\iota$ respects the restriction of $\sigma$ to
$K({\cal T})v$.

To simplify notation, we will write $F_0=K({\cal T})$.
We will now construct a valuation preserving embedding of the
henselization $F^h$ over $K$ in $(K^*,v^*)$. The restriction of this
embedding is the required embedding of $F$. Observe that $F^h$ contains
the henselization $F_0^h$. By the universal property of henselizations,
%
%
$\iota$ extends to a valuation preserving embedding of $F_0^h$ in $K^*$
since by hypothesis, $K^*$ is henselian. Since $F_0^h|F_0$ is immediate,
this embedding also respects the above mentioned restrictions of $\rho$
and $\sigma$. Through this embedding, we will from now on identify
$F_0^h$ with its image in $K^*$.

\pars
Now we have to extend $\iota$ (which by our identification has become
the identity) to an embedding of $F^h=F_0^h(a)$ in $K^*$ (over $F_0^h$)
which respects $\rho$ and $\sigma$. This is done as follows. Take a
monic polynomial $f\in {\cal O}_{F_0^h}[X]$ whose residue polynomial
$fv$ is the minimal polynomial of $av$ over $F_0^hv$; by hypothesis,
$fv$ is separable. Hensel's Lemma shows that there exists exactly one
root $a'$ of $f$ in the henselian field $K^*$ having residue
$\sigma(av)$. The assignment
\[
a\mapsto a'
\]
induces an isomorphism $\iota: F_0^h(a)\longrightarrow F_0^h(a')$ which
is valuation preserving since $F_0^h$ is henselian. As $vF^h=vF_0^h$, we
also have that $vF_0^h(a)=vF_0^h$. Thus, $\iota$ respects $\rho$ (which
after the above identification is the identity). We have to show that
$\iota$ also respects $\sigma$.

Let $n=[F_0^h(a):F_0^h]$. Since the elements $1,av,\ldots,(av)^{n-1}$ are
linearly independent, the basis $1,a,\ldots,a^{n-1}$ is a valuation
basis of $F_0^h(a)|F_0^h$, that is,
\[
v\sum_{i=0}^{n-1}c_ia^i\>=\>\min_i vc_i
\]
for any choice of $c_i\in F_0^h$. Take $g(a)\in F_0^h[a]$ where $g\in
F_0^h[X]$ is of degree $<n$; if $vg(a)=0$, then $g\in {\cal O}_{F_0^h}
[X]$ and thus, $g(a)v = (gv)(av)$. In this case,
\[
(\iota g(a))v^*=g(a')v^*= (gv)(a'v^*) =
(gv)(\sigma (av)) = \sigma((gv)(av))=\sigma(g(a)v)\;.
\]
This proves that $\iota$ respects $\sigma$.
\end{proof}

We return to the proof of Theorem~\ref{wtA}. We take any finitely
generated subextension $F|K$ of $L|K$. As pointed out above, $(F|K,v)$
is an extension without transcendence defect.
By assumption, $vK\ec vL$ and $Kv\ec Lv$, which implies that $vK\ec vF$
and $Kv\ec Fv$ because $vF|vK$ is a subextension of $vL|vK$, and $Fv|Kv$
is a subextension of $Lv|Kv$. So we can infer from Lemma~\ref{sica} that
the conditions ``$vF/vK$ is torsion free'' and ``$Fv|Kv$ is separable''
are satisfied. It now follows from Theorem~\ref{hrwtd} that $(F|K,v)$ is
strongly inertially generated. Hence by the previous lemma there is an
embedding
\[
\iota:\;(F,v) \longrightarrow (K^*,v^*)
\]
over $K$ that respects the restriction of $\rho$ to $vF$ and the
restriction of $\sigma$ to $Fv$. Since this holds for every finitely
generated subextension $(F|K,v)$ of $(L|K,v)$, it follows from part b) of
Lemma~\ref{embfingen} that also $(L,v)$ embeds in $(K^*,v^*)$ over $K$.
By Proposition~\ref{ec}, this shows that $(K,v)$ is existentially closed
in $(L,v)$, and we have now proved Theorem~\ref{wtA}.

\parm
For further use, we have to make our result more precise:

\begin{lemma} {\bf (Embedding Lemma II)}              \label{EII}\n
Take a defectless field $(K,v)$ (the valuation is allowed to be
trivial), an extension $(L|K,v)$ without transcendence defect and
an $|L|^+$-saturated henselian extension $(K^*,v^*)$ of $(K,v)$. Assume
that $vL/vK$ is torsion free and that $Lv|Kv$ is
separable. If
\[
\rho:\; vL \longrightarrow v^*K^*
\]
is an embedding over $vK$ and
\[
\sigma:\; Lv \longrightarrow K^*v^*
\]
is an embedding over $Kv$, then there exists an embedding
\[
\iota:\>(L,v) \longrightarrow (K^*,v^*)
\]
over $K$ which respects $\rho$ and $\sigma$.
\end{lemma}
\begin{proof}
Take any finitely generated subextension $(F|K,v)$ of $(L|K,v)$. Then
$(F|K,v)$ is a valued function field without transcendence defect. Since
$vL/vK$ is torsion free, the same holds for $vF/vK$. Since $Lv|Kv$ is
separable, the same holds for $Fv|Kv$. Hence by Theorem~\ref{hrwtd},
$(F|K,v)$ is strongly inertially generated, and by Embedding Lemma I,
$(F,v)$ embeds over $(K,v)$ in $(K^*,v^*)$ respecting both
embeddings $\rho$ and $\sigma$.

Using the saturation property of $(K^*,v^*)$ we wish to deduce our
assertion from this fact. To do so, we will work in an extended language
${\cal L}'$ consisting of the language ${\cal L}_{\rm VF}$ of valued
fields together with the predicates
\begin{eqnarray*}
\mbox{$\cal P$}_{\alpha}\,, & & \alpha\in \rho (vL)\\
\mbox{$\cal Q$}_{\zeta}\,,  & & \zeta\in \sigma (Lv)
\end{eqnarray*}
which are interpreted in $(K^*,v^*)$ such that
\begin{eqnarray*}
\mbox{$\cal P$}_{\alpha}(a) & \Longleftrightarrow & v^* a = \alpha\\
\mbox{$\cal Q$}_\zeta(a) & \Longleftrightarrow & av^* = \zeta
\end{eqnarray*}
for all $a\in K^*$ and in $(L,v)$ such that
\begin{eqnarray*}
\mbox{$\cal P$}_{\alpha}(b) & \Longleftrightarrow & \rho (vb) = \alpha\\
\mbox{$\cal Q$}_\zeta(b) & \Longleftrightarrow & \sigma (bv) = \zeta
\end{eqnarray*}
for all $b\in L$. Note that these interpretations coincide on $K$.
\pars
We show that $(K^*,v^*)$ remains $|L|^+$-saturated in the extended
language ${\cal L}'$. To this end, we choose a subset $S_v\subset K^*$
of representatives for all values $\alpha$ in $\rho (vL)$, and a subset
$S_r\subset K^*$ of representatives for all residues $\zeta$ in
$\sigma (Lv)$. We compute
\begin{eqnarray*}
|S_v| & = & |\rho vL| = |vL| \leq |L| < |L|^+\;,\\
|S_r| & = & |\sigma Lv| = |Lv| \leq |L| < |L|^+\;,
\end{eqnarray*}
hence $|S_v\cup S_r|<|L|^+$. Consequently, it follows
%
%
that $(K^*,v^*)$ remains $|L|^+$-saturated in the extended language
${\cal L}_{\rm VF}(S_v\cup S_r)$ (the new constants are interpreted in
$K^*$ by the corresponding elements from $S_v\cup S_r$). Now the
predicates ${\cal P}_{\alpha}$ and ${\cal Q}_\zeta$ become definable in
the language ${\cal L}_{\rm VF}(S_v\cup S_r)$. Indeed, if $\alpha\in
\rho (vL)$, then we choose $b_\alpha\in S_v$ such that $v^*b_\alpha=
\alpha$ and define ${\cal P}_\alpha(x):\Leftrightarrow v^*x=v^*
b_\alpha$. If $\zeta \in \sigma (Lv)$, then we choose $b_\zeta\in S_r$
such that $b_\zeta v^*= \zeta$ and define ${\cal Q}_\zeta(x):
\Leftrightarrow v^*(x-b_\zeta) >0$. Since $(K^*,v^*)$ is
$|L|^+$-saturated in the language ${\cal L}_{\rm VF} (S_v\cup S_r)$, it
follows that it is also $|L|^+$-saturated in the language ${\cal L}'
(S_v\cup S_r)$ and thus also in the language ${\cal L}'$, as asserted.

\pars
An embedding $\iota$ of an arbitrary subextension $(F,v)$ of
$(L|K,v)$ in $(K^*,v^*)$ over $K$ respects the predicates
${\cal P}_{\alpha}$ and ${\cal Q}_\zeta$ if and only if it
satisfies, for all $b\in F$,
\begin{eqnarray*}
&\rho (vb) = \alpha \Longleftrightarrow (F,v)\models
{\cal P}_{\alpha}(b) \Longleftrightarrow (K^*,v^*)\models
{\cal P}_{\alpha}(\iota b)\Longleftrightarrow v^*(\iota b)=\alpha\>,&\\
&\sigma (bv) = \zeta \Longleftrightarrow (F,v)\models
{\cal Q}_\zeta(b) \Longleftrightarrow (K^*,v^*)\models
{\cal Q}_\zeta(\iota b) \Longleftrightarrow (\iota b)v^* = \zeta\>,&
\end{eqnarray*}
which expresses the property of $\iota$ to respect the embeddings
$\rho$ and $\sigma$. We know that for every finitely generated
subextension of $(L|K,v)$ there exists such an embedding $\iota$. The
saturation property of $(K^*,v^*)$ now yields an embedding of $(L,v)$ in
$(K^*,v^*)$ over $K$ which respects the predicates and thus the
embeddings $\rho$ and $\sigma$. This completes the proof of our lemma.
\end{proof}

%
%
\subsection{Completions}           \label{sectcomp}
In this section, we deal with extensions of a valued field within its
completion. This is a preparation for the subsequent section on the
model theory of separably tame fields. But the results are also of
independent interest. As a preparation for the next theorem, we need:

\begin{lemma}                               \label{lCs}
Assume that $(K(x)|K,v)$ is an extension within the completion of
$(K,v)$ such that $x$ is transcendental over $K$. Then $x$ is the limit
of a pseudo Cauchy sequence in $(K,v)$ of transcendental type.
\end{lemma}
\begin{proof}
Since $x\in K^c$, it is the limit of a Cauchy sequence
$(a_\nu)_{\nu<\lambda}$ in $(K,v)$, that is, the values $v(x-a_\nu)$ are
strictly increasing with $\nu$ and are cofinal in $vK$. Suppose that
this sequence would not be of transcendental type. Then there is a
polynomial $f\in K[X]$ of least degree for which the values $vf(a_\nu)$
are not ultimately fixed. By Lemma~8 of \cite{[Ka]},
\[
vf(a_\nu)\>=\>\beta_h\,+\,hv(x-a_\nu)
\]
holds for all large enough $\nu$, where $\beta_h\in vK$ and $h$ is a
power of the characteristic exponent $p$ of $Kv$. By Lemma~9 of
\cite{[Ka]},
\[
vf(x)\> >\>\beta_h\,+\,hv(x-a_\nu)
\]
for all large enough $\nu$. As these values are cofinal in $vK$, we
conclude that $vf(x)=\infty$, that is, $f(x)=0$. Hence if $x$ is
transcendental over $K$, then $(a_\nu)_{\nu<\lambda}$ must be of
transcendental type.
\end{proof}

\begin{theorem}                                    \label{XsKc}
Let $(K,v)$ be a henselian field. Assume that $(L|K,v)$ is a separable
subextension of $(K^c|K,v)$. Then $(K,v)$ is existentially closed in
$(L,v)$. In particular, every henselian inseparably defectless field is
existentially closed in its completion.
\end{theorem}
\begin{proof}
By part a) of Lemma~\ref{embfingen}, it suffices to show that $(K,v)$ is
existentially closed in every subfield $(F,v)$ of $(L,v)$ which is
finitely generated over $K$. Equivalently, it suffices to show that
$(K,v)$ is existentially closed in $(F,v)^h$; note that $(F,v)^h\subset
(K,v)^c$ since the completion of a henselian field is again henselian
(cf.\ \cite{[W]}, Theorem~32.19).
%
%
%
As a subextension of the separable extension $L|K$, also $F|K$ is
separable. So we may choose a separating transcendence basis
${\cal T}=\{x_1,\ldots,x_n\}$ of $F|K$. Then $(F,v)$ lies in the
completion of $(K({\cal T}),v)$ since it lies in the completion of
$(K,v)$. The completion of $K({\cal T})^h$ is equal to $K^c$ since
$K({\cal T})\subseteq K^c$ and $(K^c,v)$ is henselian. Consequently,
$F^h$ lies in the completion of $K({\cal T})^h$. On the other hand,
$F^h|K({\cal T})^h$ is a finite separable extension; since a henselian
field is separable-algebraically closed in its completion (cf.\
\cite{[W]}, Theorem~32.19), it must be trivial. That is,
\[
(F,v)^h=(K(x_1,\ldots,x_n),v)^h\;.
\]

Set $F_0=K$ and $(F_i,v)=(K(x_1,\ldots,x_i),v)^h$, $1\leq i\leq n$,
where the henselization is taken within $F^h$. Now it suffices to show
that $(F_{i-1},v)\ec (F_{i},v)$ for $1\leq i\leq n$. As $x_i$ is an
element of the completion $K^c$ of $(F_{i-1},v)$, it is the limit of a
Cauchy sequence in $(F_{i-1},v)$. Since $x_i$ is transcendental
over $F_{i-1}\,$, this Cauchy sequence must be of transcendental
type by Lemma~\ref{lCs}.
%
%
Hence by Corollary~\ref{imrf}, $(F_{i-1},v)\ec (F_{i-1}(x_i),v)^h$ for
$1\leq i\leq n$, which in view of $(F_{i-1}(x_i),v)^h=(F_i,v)^h$ proves
our assertion.

The second assertion of our theorem follows from the first and the fact
that if $(K,v)$ is inseparably defectless, then the immediate extension
$K^c|K$ is separable, according to Corollary~\ref{iers}.
%
%
%
\end{proof}

From this theorem together with part b) of Lemma~\ref{sica},
we obtain:

\begin{corollary}
A henselian field $(K,v)$ is existentially closed in its completion
$K^c$ if and only if the extension $K^c|K$ is separable.
\end{corollary}
This leads to the following question:
\sn
{\bf Open Problem:} Take any field $k$. Which are the subfields
$K\subset k((t))$ with $t\in K$ such that $k((t))|K$ is separable?

\pars
Recall that $v_t$ denotes the $t$-adic valuation on $k(t)$ and on
$k((t))$. Since $(k((t)),v_t)$ is henselian, we can choose the
henselization $(k(t),v_t)^h$ in $(k((t)),v_t)$. Then $(k((t)),v_t)$ is
the completion of both $(k(t),v_t)$ and $(k(t),v_t)^h$. Further,
$(k,v_t)$ is trivially valued and thus defectless. By Theorem~\ref{ai},
it follows that $(k(t),v_t)^h$ is henselian and defectless. Now
Corollary~\ref{iers} shows:

\begin{corollary}                           \label{ciers}
The extension $k((t))|k(t)^h$ is regular.
\end{corollary}
Using Theorem~\ref{XsKc}, we conclude:
\begin{theorem}
Let $k$ be an arbitrary field. Then $(k(t),v_t)^h\ec (k((t)),v_t)$.
\end{theorem}
This result also follows from Theorem~2 of \cite{[Er6]}. It was used in
\cite{[K6]} in connection with the characterization of large fields.

\pars
To give a further application, we need another lemma.
\begin{lemma}                               \label{KecK(t)}
Let $t$ be transcendental over $K$. Suppose that $K$ admits a
nontrivial henselian valuation $v$. Then $(K,v)\ec (K(t),v_t\circ v)^h$.
\end{lemma}
\begin{proof}
Let $(K^*,v^*)$ be a $|K(t)^h|^+$-saturated elementary extension of
$(K,v)$. Then by Corollary~\ref{vgrfsat}, $v^*K^*$ is a
$|vK|^+$-saturated elementary extension of $vK$. Hence,
there exists an element $\alpha\in v^*K^*$ such that $\alpha>vK$.
%
%
We also have that $(v_t\circ v)t>vK$. Now if $\Gamma\subset \Delta$ is
an extension of ordered abelian groups and $\Delta\ni\alpha>\Gamma$,
then the ordering on $\Z\alpha+\Gamma$ is uniquely determined. Indeed,
$\Z\alpha+\Gamma$ is isomorphic to the product $\Z\alpha\amalg\Gamma$,
lexicographically ordered. So we see that the assignment $(v_t\circ
v)t\mapsto\alpha$ induces an embedding of $(v_t\circ v)K(t)\isom\Z
(v_t\circ v)t\times vK$ (with the lexicographic ordering) in $v^*K^*$
over $vK$ as ordered groups. Now choose $t^*\in K^*$ such that $v^*t^*
=\alpha$. As $(v_t\circ v)t$ and $\alpha$ are not torsion elements over
$vK$, Lemma~\ref{prelBour} shows that the assignment $t\mapsto t^*$
induces an embedding of $(K(t),v_t\circ v)$ in $(K^*,v^*)$ over $K$.
Since $(K,v)$ is henselian, so is the elementary extension $(K^*,v^*)$.
By the universal property of the henselization, the embedding can thus
be extended to an embedding of $(K(t),v_t\circ v)^h$ in $(K^*,v^*)$. By
Proposition~\ref{ec}, this gives our assertion.
\end{proof}

Now we are able to prove:
\begin{theorem}
If the field $K$ admits a nontrivial henselian valuation, then
$K\ec K((t))$ (as fields).
\end{theorem}
\begin{proof}
Let $v$ be the nontrivial valuation on $K$ for which $(K,v)$ is
henselian. By Lemma~\ref{KecK(t)}, we have that $(K,v)\ec (K(t),
v_t\circ v)^h$. By Corollary~\ref{ciers}, $K((t))|K(t)^h$ is separable.
Since $(K((t)),v_t)$ is the completion of $(K(t),v_t)$, it follows
%
%
that $(K((t)),v_t\circ v)$ is the completion of $(K(t),v_t \circ v)$.
Hence, Theorem~\ref{XsKc} shows that $(K(t),v_t\circ v)^h\ec (K((t)),
v_t\circ v)$. It follows that $(K,v)\ec (K((t)),v_t\circ v)$. In
particular, $K\ec K((t))$, as asserted.
\end{proof}

\pars
We conclude this section with the following useful result, which we will
apply in the proof of Theorem~\ref{septameAKE} in
Section~\ref{sectstsd}.

\begin{proposition}                         \label{ecdec}
Take a separable extension $(L|K,v)$ and an extension $(K_1|K,v)$ such
that $K$ is dense in $(K_1,v)$. Assume that $v$ is a valuation on
$L.K_1$ which extends the valuation $v$ from both $L$ and $K_1$ and that
$(K_1,v)\ec (L.K_1,v)$. Then $(K,v)\ec (L,v)$.
\end{proposition}
\begin{proof}
We take an $|L.K_1|^+$-saturated elementary extension
$(K_1|K,v)^*$ of the valued field extension $(K_1|K,v)$. We note that
$(K_1,v)^*$ is a subfield of the completion $K^{*c}$ of $(K,v)^*$ since
the property of $K$ to be dense in $K_1$ is elementary in the language
of valued fields with the predicate ${\cal P}$ for the subfield; indeed,
\[
\forall x \forall y \exists z:\; {\cal P}(z)\,\wedge\,
(\,y\not=0\,\rightarrow\, v(x-z)>vy\,)
\]
expresses this property.

Since $(K_1,v)\ec (L.K_1,v)$, Proposition~\ref{ec} shows that
$(L.K_1,v)$ embeds over $K_1$ in $(K_1,v)^*$. Thus $L.K_1$ can be
considered as a subfield of $K^{*c}$, and so the
same holds for the fields $L$ and $L.K^*$. Since $L|K$ is assumed to
be separable, it follows that also $L.K^*|K^*$ is separable. Now
Theorem~\ref{XsKc} shows that
\[
(K,v)^*\ec (L.K^*,v^*)\;.
\]
Since $(K,v)\prec (K,v)^*$, we obtain that $(K,v)\ec (L.K^*,v^*)$, which
yields that $(K,v)\ec (L,v)$, as asserted.
\end{proof}

%
%
\section{The Relative Embedding Property}           \label{sectREP}
Inspired by the assertion of Lemma~\ref{EII}, we define a property that
will play a key role in our approach to the model theory of tame
fields. Let {\bf C} be a class of valued fields. We will say that
{\bf C} has the \bfind{Relative Embedding Property}, if the following
holds:
\sn
if $(L,v),(K^*,v^*)\in {\bf C}$ with common subfield $(K,v)$ such that
\n
$\bullet$ \ $(K,v)$ is defectless, \n
$\bullet$ \ $(K^*,v^*)$ is $|L|^+$-saturated, \n
$\bullet$ \ $vL/vK$ is torsion free and $Lv|Kv$ is separable, \n
$\bullet$ \ there are embeddings $\rho:\; vL \longrightarrow v^*K^*$
over $vK$ and $\sigma:\; Lv \longrightarrow K^*v^*$ over $Kv$,
\n
then there exists an embedding $\iota:\>(L,v) \longrightarrow (K^*,v^*)$
over $K$ which respects $\rho$ and $\sigma$.

\pars
We will show that the Relative Embedding Property of {\bf C} implies
another property of {\bf C} which is very important for our purposes. If
${\eu C}\subset {\eu A}$ and ${\eu C}\subset {\eu B}$ are extensions of
${\cal L}$-structures, then we will write ${\eu A} \equiv_{\eu C}
{\eu B}$ if $({\eu A},{\eu C}) \equiv ({\eu B},{\eu C})$ in the language
${\cal L}({\eu C})$ augmented by constant names for the elements of
${\eu C}$. If for every two fields $(L,v),(F,v)\in{\bf C}$ and every
common defectless subfield $(K,v)$ of $(L,v)$ and $(F,v)$ such
that $vL/vK$ is torsion free and $Lv|Kv$ is separable, the side
conditions $vL\equiv_{vK} vF$ and $Lv\equiv_{Kv} Fv$ imply that
$(L,v)\equiv_{(K,v)} (F,v)$, then we will call {\bf C} \bfind{relatively
subcomplete}. Note that if {\bf C} is a relatively subcomplete class of
defectless fields, then {\bf C} is relatively model complete: the side
conditions $vK\prec vL$ and $Kv\prec Lv$ imply that $vL/vK$ is torsion
free and $Lv|Kv$ is separable (by Lemma~\ref{sica}) and that $vK
\equiv_{vK} vL$ and $Kv\equiv_{Kv}Lv$, hence if {\bf C} is relatively
subcomplete, then we obtain $(K,v)\equiv_{(K,v)} (L,v)$, that is,
$(K,v)\prec (L,v)$. But relative model completeness is weaker than
relative subcompleteness, because $vL\equiv_{vK} vF$ does not imply that
$vK\prec vL$, and $Lv\equiv_{Kv} Fv$ does not imply that $Kv\prec Lv$.

The following lemma shows that the Relative Embedding Property is a
powerful property:
\begin{lemma}                               \label{REPimp}
Take an elementary class {\bf C} of defectless valued fields which has
the Relative Embedding Property. Then {\bf C} is relatively subcomplete
and relatively model complete, and the AKE$^\exists$ Principle is
satisfied by all extensions $(L|K,v)$ such that both $(K,v),(L,v)\in
{\bf C}$. If moreover all fields in {\bf C} are of fixed equal
characteristic, then {\bf C} is relatively complete.
\end{lemma}
\begin{proof}
%
Let us first show that $(L|K,v)$ satisfies the AKE$^\exists$ Principle
whenever $(K,v),(L,v)\in {\bf C}$. So assume that $vK\ec vL$ and $Kv\ec
Lv$. We take an $|L|^+$-saturated elementary extension $(K^*,v^*)$
of $(K,v)$. Since {\bf C} is assumed to be an elementary class,
$(K,v)\in {\bf C}$ implies that $(K^*,v^*)\in {\bf C}$. Because of
$vK\ec vL$ and $Kv\ec Lv$, there are embeddings $vL\,\rightarrow\,
v^*K^*$ over $vK$ and $Lv\,\rightarrow\,K^*v^*$ over $Kv$ by
Proposition~\ref{ec}. Moreover, $vL/vK$ is torsion free and $Lv|Kv$ is
separable by Lemma~\ref{sica}. So by the Relative Embedding Property
there is an embedding of $(L,v)$ in $(K^*,v^*)$ over $K$, which shows
that $(K,v)\ec (L,v)$.

\pars
In order to show that {\bf C} is relatively subcomplete, we take $(L,v),
(F,v)\in {\bf C}$ with common defectless subfield $(K,v)$ such that
$vL/vK$ is torsion free, $Lv|Kv$ is separable, $vL\equiv_{vK} vF$ and
$Lv\equiv_{Kv} Fv$. We have to show that $(L,v)\equiv_{(K,v)} (F,v)$.

To begin with, we construct an elementary extension $(L_0,v)$ of $(L,v)$
and an elementary extension $(F_0,v)$ of $(F,v)$ such that $vL_0=vF_0$
and $L_0v=F_0v$. Our condition $vL\equiv_{vK} vF$ means that $vL$ and
$vF$ are equivalent in the augmented language ${\cal L}_{\rm OG}(vK)$ of
ordered groups with constants from $vK$. Similarly, $Lv\equiv_{Kv} Fv$
means that $Lv$ and $Fv$ are equivalent in the augmented language ${\cal
L}_{\rm R}$ of rings with constants from $Kv$. It follows from the proof
of Theorem 6.1.15 in \cite{[C--K]} that we can
choose a cardinal $\lambda$ and an ultrafilter ${\cal D}$ on $\lambda$
such that $\prod_{\lambda} vL/{\cal D}\isom \prod_{\lambda} vF/{\cal D}$
and $\prod_{\lambda} Lv/{\cal D}\isom \prod_{\lambda} Fv/{\cal D}$ in
the respective augmented languages. But this means that for $(L_0,v):=
\prod_{\lambda} (L,v)/{\cal D}$ and $(F_0,v):= \prod_{\lambda}
(F,v)/{\cal D}$, we have that $vL_0=\prod_{\lambda} vL/{\cal D}$ is
isomorphic over $vK$ to $vF_0=\prod_{\lambda} vF/{\cal D}$, and
$L_0v=\prod_{\lambda} Lv/{\cal D}$ is isomorphic over $Kv$ to
$F_0v=\prod_{\lambda} vF/{\cal D}$. Passing to an equivalent valuation
on $L_0$ which still extends the valuation $v$ of $K$, we may assume
that $vL_0=vF_0$; similarly, passing to an equivalent residue map we
may assume that $L_0v=F_0v$. As $vL/vK$
%
%
is torsion free by assumption and $vL_0/vL$
%
%
are torsion free since $vL\prec vL_0\,$,
%
%
we find that $vL_0/vK=vF_0/vK$ is torsion free. Similarly, one shows
that $L_0v=F_0v$ is a separable extension of $Kv$.

\pars
Now we construct two elementary chains $((L_i,v))_{i<\omega}$ and
$((F_i,v))_{i<\omega}$ as follows. We choose a cardinal $\kappa_0
=\max\{|L_0|,|F_0|\}$. By induction, for every $i<\omega$ we take
$(L_{i+1},v)$ to be a $\kappa_i^+$-saturated elementary extension of
$(L_i,v)$, where $\kappa_i =\max\{|L_i|,|F_i|\}$, and $(F_{i+1},v)$ to
be a $\kappa_i^+$-saturated elementary extension of $(F_i,v)$. We can
take $(L_{i+1},v)=\prod_{\lambda_i} (L_i,v)/{\cal D}_i$ and
$(F_{i+1},v)=\prod_{\lambda_i} (F_i,v)/{\cal D}_i$ for suitable
cardinals $\lambda_i$ and ultrafilters ${\cal D}_i$; this yields that
$vL_i=vF_i$ and $L_iv=F_iv$ for all $i$.

All $(L_i,v)$ and $(F_i,v)$ are elementary extensions of $(L,v)$ and
$(F,v)$ respectively, so it follows that they lie in {\bf C} and in
particular, are defectless fields. We take $(L^*,v)$ to be the union
over the elementary chain $(L_i,v)$, $i<\omega$; so $(L,v)\prec
(L^*,v)$. Similarly, we take $(F^*,v)$ to be the union over the
elementary chain $(F_i,v)$, $i<\omega$; so $(F,v)\prec (F^*,v)$. Now
we carry out a back and forth construction that will show that $(L^*,v)$
and $(F^*,v)$ are isomorphic over $K$.

We start by embedding $(L_0,v)$ in $(F_1,v)$. The identity mappings are
embeddings of $vL_0$ in $vF_1$ over $vK$ and of $L_0v$ in $F_1v$ over
$Kv$, and we know that $vL_0/vK$ is torsion free and $L_0v|Kv$ is
separable. Since $(F_1,v)$ is $\kappa_0^+$-saturated with $\kappa_0\geq
|L_0|$, and since $(K,v)$ is defectless, we can apply the Relative
Embedding Property to find an embedding $\iota_0$ of $(L_0,v)$ in
$(F_1,v)$ over $K$ which respects the embeddings of the value group and
the residue field. That is, we have that $v\iota_0 L_0 = vF_0$ and
$(\iota_0 L_0)v= F_0 v\,$.

The isomorphism $\iota_0^{-1}:\>\iota_0 L_0\,\rightarrow\,L_0$ can be
extended to an isomorphism $\iota_0^{-1}$ from $F_1$ onto some
extension field of $L_0$ which we will simply denote by $\iota_0^{-1}
F_1\,$. We take the valuation on this field to be the one induced via
$\iota_0^{-1}$ by the valuation on $F_1$. Hence, $\iota_0^{-1}$ induces
an isomorphism on the value groups and the residue fields, so that we
obtain that $v\iota_0^{-1}F_1=vF_1=vL_1$ and $(\iota_0^{-1} F_1)v=F_1v=
L_1v$. The identity mappings are embeddings of $v\iota_0^{-1}F_1$ in
$vL_2$ over $vL_0$ and of $(\iota_0^{-1}F_1)v$ in $L_2v$ over $L_0v$.
Since $vL_0\prec vL_1$ and $L_0v\prec L_1v$, we know that $v\iota_0^{-1}
F_1/vL_0$ is torsion free and $(\iota_0^{-1} F_1)v|L_0v$ is separable.
Since $(L_2,v)$ is $\kappa_1^+$-saturated with $\kappa_1\geq |F_1|=
|\iota_0^{-1} F_1|$, and since $(L_0,v)$ is defectless, we can apply the
Relative Embedding Property to find an embedding $\tilde{\iota}_1$ of
$(\iota_0^{-1} F_1,v)$ in $(L_2,v)$ over $L_0$ which respects the
embeddings of the value group and the residue field. That is, we obtain
an embedding $\iota'_1:=\tilde{\iota}_1\iota_0^{-1}$ of $F_1$ in $L_2$
over $K$.
%
%
We note that $L_0\subset \iota'_1 F_1$ and that ${\iota'_1}^{-1}:
\iota'_1 F_1\rightarrow F_1$ extends $\iota_0\,$.

\pars
Suppose that we have constructed, for an even $i$, the embeddings
\begin{eqnarray*}
\iota_i: & & (L_i,v)\;\longrightarrow\;(F_{i+1},v)\\
{\iota'_{i+1}}: & & (F_{i+1},v) \;\longrightarrow\; (L_{i+2},v)
\end{eqnarray*}
as embeddings over $K$, such that $L_i\subset {\iota'_{i+1}} F_{i+1}$ and
that ${\iota'_{i+1}}^{-1}: \iota'_{i+1} F_{i+1}\rightarrow F_{i+1}$
extends $\iota_i\,$. We wish to construct similar embeddings for $i+2$
in place of $i$.

The isomorphism ${\iota'_{i+1}}^{-1}: \iota'_{i+1} F_{i+1}\rightarrow
F_{i+1}$ can be extended to an isomorphism ${\iota'_{i+1}}^{-1}$ from
$L_{i+2}$ onto some extension field of $F_{i+1}$ which we will denote by
${\iota'_{i+1}}^{-1} L_{i+2}\,$; this isomorphism extends $\iota_i\,$.
We take the valuation on this field to be the one induced via
${\iota'_{i+1}}^{-1}$ by the valuation on $L_{i+2}$. We obtain that
$v{\iota'_{i+1}}^{-1} L_{i+2}= vL_{i+2}=vF_{i+2}$ and
$({\iota'_{i+1}}^{-1} L_{i+2})v=L_{i+2}v= F_{i+2}v$. The identity
mappings are embeddings of $v{\iota'_{i+1}}^{-1} L_{i+2}$ in $vF_{i+3}$
over $vF_{i+1}$ and of $({\iota'_{i+1}}^{-1} L_{i+2})v$ in $F_{i+3}v$
over $F_{i+1}v$. Since $vF_{i+1}\prec vF_{i+3}$ and $F_{i+1}v\prec
F_{i+3}v$, we know that $v{\iota'_{i+1}}^{-1} L_{i+2} /vF_{i+1}$ is
torsion free and $({\iota'_{i+1}}^{-1} L_{i+2})v|F_{i+1}v$ is separable.
Since $(F_{i+3},v)$ is $\kappa_{i+2}^+$-saturated with $\kappa_{i+2}\geq
|L_{i+2}|= |{\iota'_{i+1}}^{-1} L_{i+2}|$, and since $(F_{i+1},v)$ is
defectless, we can apply the Relative Embedding Property to find an
embedding $\tilde{\iota}'_{i+2}$ of $({\iota'_{i+1}}^{-1} L_{i+2},v)$ in
$(F_{i+3},v)$ over $F_{i+1}$ which respects the embeddings of the value
group and the residue field. We obtain an embedding $\iota_{i+2}:=
\tilde{\iota}'_{i+2} {\iota'_{i+1}}^{-1}$ of $L_{i+2}$ in $F_{i+3}$;
since $\tilde{\iota}'_{i+2}$ is the identity on $\iota_i L_i\subset
F_{i+1}$ and ${\iota'_{i+1}}^{-1}$ extends $\iota_i$, this embedding
also extends $\iota_i$.
%
%
We note that $F_{i+1}\subset \iota_{i+2} L_{i+2}$ and that
$\iota_{i+2}^{-1}: \iota_{i+2} L_{i+2}\rightarrow L_{i+2}$
extends ${\iota'_{i+1}}\,$.
%

\pars
Now we take $\iota$ to be the set theoretical union over the embeddings
$\iota_i\,$, $i<\omega$ even. Then $\iota$ is an embedding of $(L^*,v)$
in $(F^*,v)$. It is onto since $F_i$ lies in the image of
$\iota_{i+1}$, for every odd $i$. So we have obtained an isomorphism
from $(L^*,v)$ onto $(F^*,v)$ over $K$, which shows that $(L^*,v)
\equiv_{(K,v)} (F^*,v)$. Since $(L,v)\prec (L^*,v)$ and $(F,v)\prec
(F^*,v)$, this implies that $(L,v) \equiv_{(K,v)} (F,v)$, as required.
We have proved that {\bf C} is relatively subcomplete, and we know
already that this implies that {\bf C} is relatively model complete.

\pars
Finally, assume in addition that all fields in {\bf C} are of fixed
equal characteristic. We wish to show that {\bf C} is relatively
complete. So take $(L,v),(F,v)\in {\bf C}$ such that $vL\equiv vF$ and
$Lv\equiv Fv$. Fixed characteristic means that $L$ and $F$ have
a common prime field $K$. The assumption that both $(L,v)$ and $(F,v)$
are of equal characteristic means that the restrictions of their
valuations to $K$ is trivial. Hence, $vK=0$ and consequently, $vL/vK$ is
torsion free and $vL\equiv vF$ implies that $vL\equiv_{vK} vF$. Further,
$K=Kv$ is also the prime field of $Lv$ and $ Fv$, so
$Lv\equiv Fv$ implies that $Lv\equiv_{Kv} Fv$.
%
%
Since a prime field is always perfect, we also have that $Lv|Kv$ is
separable. As a trivially valued field, $(K,v)$ is defectless. From what
we have already proved, we obtain that $(L,v)\equiv_{(K,v)} (F,v)$,
which implies that $(L,v)\equiv (F,v)$.
\end{proof}

Now we look for a criterion for an elementary class of valued fields to
have the Relative Embedding Property. In some way, we have to improve
Embedding Lemma II (Lemma~\ref{EII}) to cover the case of extensions
$(L|K,v)$ with transcendence defect. Loosely speaking, these contain an
immediate part. The idea is to require that this part can be treated
separately, that is, that we find an intermediate field $(L',v)\in
{\bf C}$ such that $(L|L',v)$ is immediate and $(L'|K,v)$ has no
transcendence defect. The immediate part has then to be handled by a new
approach which we will describe in the following embedding lemma. Note
that by Theorem~1 of \cite{[Ka]} together with Theorem~\ref{fix}, the
hypothesis on $x$ does automatically hold if $(K,v)$ is algebraically
maximal.

\begin{lemma} {\bf (Embedding Lemma III)} \label{EIII} \n
Let $(K(x)|K,v)$ be a nontrivial immediate extension of valued
fields. If $x$ is the limit of a pseudo Cauchy sequence of
transcendental type in $(K,v)$,
%
%
then $(K(x),v)^h$ embeds over $K$ in every $|K|^+$-saturated
henselian extension $(K,v)^{*}$ of $(K,v)$.
\end{lemma}
\begin{proof}
Take a pseudo Cauchy sequence $(a_\nu)_{\nu<\lambda}$ of transcendental
type in $(K,v)$ with limit $x$. Then the collection of elementary
formulas ``$v(x-a_\nu)=v(a_{\nu+1}-a_\nu)$'', $\nu<\lambda$, is a
(partial) type over $(K,v)$. Indeed, if a finite subset of these
formulas is given and $\nu_0$ is the largest of the indeces $\nu$, then
all formulas in the subset are satisfied by $x=a_{\nu+1}$.

Since $(K,v)^{*}$ is $|K|^+$-saturated, there is an element $x^{*}\in
K^{*}$ such that $v^{*}(x^{*}-a_\nu)= v^{*}(a_{\nu+1}-a_\nu)$ holds for
all $\nu<\lambda$. That is, $x^{*}$ is also a limit of
$(a_\nu)_{\nu<\lambda}$. By Theorem~2 of \cite{[Ka]}, the homomorphism
induced by $x \mapsto x^{*}$ is an embedding of $(K(x),v)$ over $K$ in
$(K,v)^{*}$. By the universal property of the henselization,
%
%
this embedding can be extended to an embedding of $(K(x),v)^h$ over
$K$ in $(K,v)^{*}$, since the latter is henselian by hypothesis.
\end{proof}

\n
Note that the lemma fails if the condition on the
%
%
pseudo Cauchy sequence to be transcendental is omitted, even if we
require in addition that $(K,v)$ is henselian.
%
%
There may exist nontrivial finite immediate extensions $(K(x)|K,v)$ of
henselian fields; for a comprehensive collection of examples, see
\cite{[K9]}. On the other hand, $K^*$ may be a regular extension of $K$
(e.g., this is always the case if $(K,v)^*$ is an elementary extension
of $(K,v)\,$), and then, $K(x)$ does certainly not admit an embedding
over $K$ in $K^*$.

\parm
The model theoretic application of Embedding Lemma~III is:
\begin{corollary}                          \label{imrf}
Let $(K,v)$ be a henselian field and $(K(x)|K,v)$ an immediate extension
such that $x$ is the limit of a pseudo Cauchy sequence of transcendental
type in $(K,v)$. Then $(K,v)\ec (K(x),v)^h$. In particular, an
algebraically maximal field is existentially closed in every
henselization of an immediate rational function field of
transcendence degree $1$.
\end{corollary}
\begin{proof}
Choose $(K,v)^*$ to be a $|K|^+$-saturated elementary extension of
$(K,v)$. Since ``henselian'' is an elementary property, $(K,v)^*$ will
also be henselian. Now apply Embedding Lemma III and
Proposition~\ref{ec}.
\end{proof}

Now we are able to give the announced criterion:
\begin{lemma}                               \label{critREP}
Let {\bf C} be an elementary class of valued fields which satisfies
\begin{axiom}
\ax{(CALM)} every field in {\bf C} is algebraically maximal,
\ax{(CRAC)} if $(L,v)\in {\bf C}$ and $K$ is relatively algebraically
closed in $L$ such that $Lv| Kv$ is algebraic and $vL/vK$ is
a torsion group, then $(K,v)\in {\bf C}$ with $Lv= Kv$ and
$vL=vK$,
\ax{(CIMM)} if $(K,v)\in {\bf C}$, then every henselization of an
immediate function field of transcendence degree~$1$ over $(K,v)$ is
already the henselization of a rational function field over $K$.
\end{axiom}
Then {\bf C} has the Relative Embedding Property.
\end{lemma}
\begin{proof}
Assume that the elementary class {\bf C} satisfies (CALM), (CRAC) and
(CIMM). Take $(L,v), (K^*,v^*)\in {\bf C}$ with $(K^*,v^*)$ being
$|L|^+$-saturated, a defectless valued subfield $(K,v)$ of $(L,v)$ and
$(K^*,v^*)$ such that $vL/vK$ is torsion free and $Lv|Kv$ is separable,
and embeddings $\rho:\; vL \rightarrow v^*K^*$ over $vK$ and $\sigma:\>
Lv \rightarrow K^*v^*$ over $Kv$. We have to show that there exists an
embedding $\iota:\>(L,v) \rightarrow (K^*,v^*)$ over $K$ which respects
$\rho$ and $\sigma$.

Take the set ${\cal T}=\{x_i\,,\,y_j\mid i\in I\,,\,j\in J\}$ as in the
proof of Corollary~\ref{tint}. Then $vL/vK({\cal T})$ is a torsion group
and $Lv| K({\cal T})v$ is algebraic. Let $K'$ be the relative algebraic
closure of $K({\cal T})$ within $L$. It follows that also $vL/vK'$ is a
torsion group and $Lv| K'v$ is algebraic. Hence by condition (CRAC), we
have that $(K',v)\in {\bf C}$ with $Lv= K'v$ and $vL=vK'$, which shows
that the extension $L|K'$ is immediate. On the other hand, ${\cal T}$ is
a standard valuation transcendence basis of $(K'|K,v)$ by construction,
hence according to Corollary~\ref{svtb->wtd}, this extension has no
transcendence defect. Since $(K,v)$ is defectless by assumption and
$(K^*,v^*)$ is henselian by condition (CALM), Lemma~\ref{EII} gives an
embedding of $(K',v)$ in $(K^*,v^*)$ over $K$ which respects $\rho$ and
$\sigma$. Now we have to look for an extension of this embedding to
$(L,v)$. Since $(L|K',v)$ is immediate, such an extension will
automatically respect $\rho$ and $\sigma$.

We identify $K'$ with its image in $K^*$. In view of part b) of
Lemma~\ref{embfingen}, it remains to show that every finitely generated
subextension $(F,v)$ of $(L|K',v)$ embeds over $K'$ in $(K^*,v^*)$.
We apply our slicing approach. Since $F$ is finitely generated
over $K'$, it has a finite transcendence basis $\{t_1,\ldots,t_n\}$ over
$K'$. Let us put $K_0=K'$ and $K_i$ to be the relative algebraic closure
of $K(t_1, \ldots,t_i)$ in $L$ for $1\leq i\leq n$. Then $K_n$ contains
$F$, and by condition (CRAC), every $(K_i,v)$ is a member of {\bf C}.
Moreover, $\mbox{\rm trdeg}(K_{i+1}|K_i)=1$ for $0\leq i<n$. We proceed
by induction on $i$. If we have shown that $(K_i,v)$ embeds in
$(K^*,v^*)$ over $K'$, then we identify it with its image. Hence it
now remains to show that the immediate extension $(K_{i+1},v)$ of
transcendence degree 1 embeds in $(K^*,v^*)$ over $K_i\,$. Since
$(K^*,v^*)$ is $|L|^+$-saturated, it is also $|K_{i+1}|^+$-saturated.
Hence again, part b) of Lemma~\ref{embfingen} shows that it suffices to
prove the existence of an embedding for every finitely generated
subextension $(F_{i+1},v)$ of $(K_{i+1}| K_i,v)$. Since $(F_{i+1}|
K_i,v)$ is an
immediate function field of transcendence degree 1, by condition (CIMM),
its henselization is the henselization $K_i(x_{i+1})^h$ of a rational
function field. Since $(K_i,v)$ is algebraically maximal by condition
(CALM), Theorem~\ref{fix} shows that $x_{i+1}$ is the limit of a pseudo
Cauchy sequence of transcendental type in $(K_i,v)$. Now Embedding Lemma
III (Lemma~\ref{EIII}) now yields that there is an embedding of
$(F_{i+1},v)$ in $(K^*,v^*)$ over $K_i\,$. This completes our proof by
induction.
\end{proof}

%
%
\section{The model theory of tame and separably tame fields}
\label{sectmtta}
%
%
\subsection{Tame fields}                    \label{sectmttf}
We have already shown in part a) of Corollary~\ref{cortame} that in
positive characteristic, the class of tame fields coincides with the
class of algebraically maximal perfect fields. Let us show that the
property of being a tame field of fixed residue characteristic is
elementary. If the residue characteristic is fixed to be $0$ then by
Theorem~\ref{tamerc0}, ``tame'' is equivalent to ``henselian'' which is
axiomatized by the axiom scheme (HENS). Now assume that the residue
characteristic is fixed to be a positive prime $p$. By
Theorem~\ref{tame}, a valued field of
positive residue characteristic is tame if and only if it is an
algebraically maximal field having $p$-divisible value group and perfect
residue field. A valued field $(K,v)$ has $p$-divisible value group if
and only if it satisfies the following elementary axiom:
\begin{axiom}
\ax{(VGD$_p$)} $\forall x\,\exists y:\> vxy^p = 0\>\vee\>x=0\;.$
\end{axiom}
Furthermore, $(K,v)$ has perfect residue field if and only if it
satisfies:
\begin{axiom}
\ax{(RFD$_p$)} $\forall x\,\exists y:\> vx=0\>\rightarrow\>v(xy^p - 1) >
0\;.$
\end{axiom}
Finally, the property of being algebraically maximal is axiomatized by
the axiom schemes (HENS) and (MAXP). We summarize: The \bfind{theory of
tame fields of residue characteristic $0$} is just the theory of
henselian fields of residue characteristic $0$. If $p$ is a prime, then
the \bfind{theory of tame fields of residue characteristic $p$} is the
theory of valued fields together with axioms (VGD$_p$), (RFD$_p$),
(HENS) and (MAXP). Now we also see how to axiomatize the theory of all
tame fields. Indeed, for residue characteristic 0 there are no
conditions on the value group and the residue field. For residue
characteristic $p>0$, we have to require (VGD$_p$) and (RFD$_p$). We can
do this by the axiom scheme
\begin{axiom}
\ax{(TAD)} $v(\underbrace{1+\ldots+1}_{p \ \rm times})>0\>
\rightarrow\>\mbox{\rm (VGD$_p$)}\,\wedge\,
\mbox{\rm (RFD$_p$)}\hfill(p\mbox{\rm\ prime})\,. \hspace{2cm}$
\end{axiom}
So the \bfind{theory of tame fields} is the theory of valued fields
together with axioms (TAD), (HENS) and (MAXP).

\pars
Recall that by part a) of Corollary~\ref{cortame}, a valued field of
positive characteristic is tame if and only if it is algebraically
maximal and perfect. We have already seen in Lemma~\ref{AKEex-am} that
every AKE$^\exists$-field must be algebraically maximal. Therefore, the
model theory of tame fields that we will develop now is representative
of the model theory of perfect valued fields in positive characteristic.

Let {\bf C} be the elementary class of all tame fields. By
Lemma~\ref{thdp}, all tame fields are henselian defectless, so {\bf C}
satisfies condition (CALM) of Lemma~\ref{critREP}. By Lemma~\ref{trac},
it also satisfies condition (CRAC). Finally, it satisfies (CIMM) by
virtue of Theorem~\ref{stt3}. Hence, we can infer from
Lemma~\ref{critREP} and Lemma~\ref{REPimp}:

\begin{theorem}                             \label{tamemod1}
The elementary class of tame fields has the Relative Embedding Property
and is relatively subcomplete and relatively model complete.
Every elementary class of tame fields of fixed equal characteristic
is relatively complete.
\end{theorem}

Lemma~\ref{critREP} does not give the full information about the
AKE$^\exists$  Principle because it requires that not only $(K,v)$, but
also $(L,v)$ is a member of the class {\bf C}. If the latter is not the
case, then it just suffices if one can show that it is contained in a
member of {\bf C}. To this end, we need the following lemma:
%
%
\begin{lemma}           \label{te}
If $\Gamma$ is a $p$-divisible ordered abelian group and $\Gamma\ec
\Delta$, then $\Gamma$ is also existentially closed in the
$p$-divisible hull of $\Delta$. If $k$ is a perfect field and $k\ec
\ell$, then $k$ is also existentially closed in the perfect hull of
$\ell$.

If $(K,v)$ is a tame field and $(L|K,v)$ an extension with $vK\ec
vL$ and $Kv\ec Lv$, then every maximal purely wild extension $(W,v)$ of
$(L,v)$ is a tame field satisfying $vK\ec vW$ and $Kv\ec Wv$.
\end{lemma}
\begin{proof}
By Proposition~\ref{ec}, $\Gamma\ec \Delta$ implies that $\Delta$
embeds over $\Gamma$ in every $|\Delta|^+$-saturated elementary
extension of $\Gamma$. Such an elementary extension is $p$-divisible
like $\Gamma$. Hence, the embedding can be extended to an embedding of
$\frac{1}{p^{\infty}}\Delta$, which by Proposition~\ref{ec} shows that
$\Gamma\ec\frac{1}{p^{\infty}}\Delta$.

Again by the same lemma, $k\ec\ell$ implies that $\ell$ embeds
over $k$ in every $|\ell|^+$-saturated elementary extension of $k$. Such
an elementary extension is perfect like $k$. Hence, the embedding can be
extended to an embedding of $\ell^{1/p^{\infty}}$, which by
Proposition~\ref{ec} shows that $k\ec\ell^{1/p^{\infty}}$.

Now suppose that the assumptions of the final assertion of our
lemma hold. By Corollary~\ref{Wtame}, $(W,v)$ is a tame field. By
Theorem~\ref{P}, $vW$ is the $p$-divisible hull $\frac{1}{p^{\infty}}vL$
of $vL$, and $ Wv$ is the perfect hull $Lv^{1/p^{\infty}}$ of $Lv$. So
our assertion follows since we have just proved that $vK$ (which is
$p$-divisible by Theorem~\ref{tame}) is existentially closed in
$\frac{1}{p^{\infty}}vL$ and that $Kv$ (which is perfect by
Theorem~\ref{tame}) is existentially closed in the perfect hull
$Lv^{1/p^{\infty}}$ of $Lv$.
\end{proof}

Assume that $(K,v)$ is a tame field and $(L|K,v)$ an extension such that
$vK\ec vL$ and $Kv\ec Lv$. We choose some maximal purely wild extension
$(W,v)$ of $(L,v)$. According to the foregoing lemma, $(W,v)$ is a tame
field with $vK\ec vW$ and $Kv\ec Wv$. Hence by Theorem~\ref{tamemod1}
together with Lemma~\ref{REPimp}, $(K,v)\ec (W,v)$. It follows that
$(K,v)\ec (L,v)$. This proves the first assertion of
Theorem~\ref{tameAKE}.

\pars
Now let {\bf C} be an elementary class of valued fields. We define
\[v{\bf C}\,:=\,\{vK\mid (K,v)\in {\bf C}\}\;\mbox{\ \ and\ \ }\;
{\bf C}v\,:=\,\{Kv\mid (K,v)\in {\bf C}\}\;.\]
If both $v{\bf C}$ and ${\bf C}v$ are model complete elementary
classes, then the side conditions $vK\prec vL$ and $Kv\prec Lv$
will hold for every two members $(K,v)\subset (L,v)$ of {\bf C}.
Similarly, if $v{\bf C}$ and ${\bf C}v$ are complete elementary
classes, then the side conditions $vK\equiv vL$ and $Kv\equiv
Lv$ will hold for all $(K,v),(L,v)\in {\bf C}$. So we obtain from
the foregoing theorems:
\begin{theorem}                           \label{cormodtf}
If {\bf C} is an elementary class consisting of tame fields and if
$v{\bf C}$ and ${\bf C}v$ are model complete elementary classes, then
{\bf C} is model complete. If {\bf C} is an elementary class consisting
of tame fields of fixed equal characteristic, and if $v{\bf C}$ and
${\bf C}v$ are complete elementary classes, then {\bf C} is complete.
\end{theorem}
Note that the converses are true by virtue of Corollary~\ref{vgrfequiv},
provided that $v{\bf C}$ and ${\bf C}v$ are elementary classes.

\pars
\begin{corollary}
Let {\bf T} be an elementary theory consisting of all perfect valued
fields of equal characteristic whose value groups satisfy a given model
complete elementary theory ${\bf T}_{\rm vg}$ of ordered abelian groups
and whose residue fields satisfy a given model complete elementary
theory ${\bf T}_{\rm rf}$ of fields. Then the theory {\bf T}$^*$ of
algebraically maximal valued fields satisfying {\bf T} is the model
companion of {\bf T}.
\end{corollary}
\begin{proof}
It follows from Theorem~\ref{cormodtf} and Corollary~\ref{cortame} that
{\bf T}$^*$ is model complete. For every model $K$ of {\bf T}, any
maximal immediate algebraic extension is a model of {\bf T}$^*$ because
it has the same value group and residue field.
\end{proof}
\n
%
In the case of positive characteristic, ${\bf T}^*$ is in general not a
model completion since there exist perfect valued fields of positive
characteristic which admit two nonisomorphic maximal immediate algebraic
extensions, both being models of the model companion. In the case of
equal characteristic 0, the algebraically maximal fields are just the
henselian fields, and
%
%
we find that ${\bf T}^*$ is a model completion of {\bf T}, because
henselizations are unique up to isomorphism.

\parm
A \bfind{weak prime model} in an elementary class {\bf C} is a model in
{\bf C} that can be embedded in every other highly enough saturated
member of {\bf C}. Elementary classes of tame fields of equal
characteristic admit weak prime models if the elementary classes of
their value groups and their residue fields do:

\begin{theorem}                     \label{prmo}
Let {\bf C} be an elementary class consisting of tame fields of
equal characteristic. Suppose that there exists an infinite cardinal
$\kappa$, an ordered group $\Gamma$ and a field $k$, both of cardinality
$\leq\kappa$, such that $\Gamma$ admits an elementary embedding in every
$\kappa^+$-saturated member of $v{\bf C}$ and $k$ admits an elementary
embedding in every $\kappa^+$-saturated member of ${\bf C}v$. Then there
exists $(K_0,v)\in {\bf C}$ of cardinality $\leq\kappa$, having value
group $\Gamma$ and residue field $k$, such that $(K_0,v)$ admits an
elementary embedding in every $\kappa^+$-saturated member of {\bf C}.
Moreover, we can assume that $(K_0,v)$ admits a standard valuation
transcendence basis over its prime field.
\end{theorem}
\begin{proof}
%
Take any $(E,v)\in {\bf C}$ and let $(E,v)^*$ be a $\kappa^+$-saturated
elementary extension of $(E,v)$. Then also $v^*E^*$ and $E^*v^*$ are
$\kappa^+$-saturated. Since {\bf C} is an elementary class, we find that
$(E,v)^*\in {\bf C}$. Consequently, $(E,v)^*$ is a tame field. By
Theorem~\ref{tame}, its value group is $p$-divisible and its residue field
is perfect. By assumption, $\Gamma$ admits an elementary embedding in
$v^*E^*$, and $k$ admits an elementary embedding in $E^*v^*$. Hence,
also $\Gamma$ is $p$-divisible and $k$ is perfect.

Now by Lemma~\ref{Gk}, there exists a tame field $(K_0,v)$ of the same
characteristic as $k$ and cardinality at most $\kappa$, having value
group $\Gamma$ and residue field $k$ and admitting a standard valuation
transcendence basis over its prime field. If $(K^*,v^*)$ is a
$\kappa^+$-saturated model of {\bf C}, then $v^*K^*$ and $K^*v^*$ are
$\kappa^+$-saturated models of $v{\bf C}$ and ${\bf C}v$ respectively.
Hence by hypothesis, there exists an elementary embedding of $\Gamma$ in
$v^*K^*$ over the trivial group $\{0\}$, and an elementary embedding of
$k$ in $K^*v^*$ over the prime field $k_0$ of $k$. Now $k_0$ is at the
same time the prime field of $K_0v$ and of $K^*v^*$. As we are dealing
with valued fields of equal characteristic, $k_0$ is also the prime
field of $K_0$ and $K^*$, and the valuation $v$ is trivial on $k_0\,$.
We have that $vK_0/vk_0$ is torsion free and $K_0v|k_0v$ is separable.
Now Embedding Lemma II (Lemma~\ref{EII}) shows the existence of an
embedding of $(K_0,v)$ in $(K^*,v^*)$ over $k_0\,$. By virtue of
Theorem~\ref{tamemod1}, this embedding is elementary (because the
embeddings of value group and residue field are). This shows that
$(K_0,v)$ is elementarily embeddable in every $\kappa^+$-saturated model
of {\bf C}. This in turn shows that $(K_0,v)$ is a model of {\bf C} and
thus a weak prime model of {\bf C}.
\end{proof}

\pars
The weak prime models that we have constructed in the foregoing proof
have the special property that they admit a standard valuation
transcendence basis over their prime field. The following corollary
confirms the representative role of models with this property.

\begin{corollary}
For every tame field $(L,v)$ of arbitrary characteristic, there exists a
tame subfield $(K,v)\prec (L,v)$ such that $(K,v)$ admits a standard
valuation transcendence basis over its prime field and $(L|K,v)$ is
immediate.
\end{corollary}
\begin{proof}
According to Corollary~\ref{tint}, for every tame field $(L,v)$ there
exists a tame subfield $(K,v)$ of $(L,v)$ admitting a standard valuation
transcendence basis over its prime field, such that $(L|K,v)$ is
immediate. In view of Theorem~\ref{tamemod1}, the latter fact shows that
$(K,v)\prec (L,v)$.
\end{proof}

\pars
As an example, we consider the theory of tame fields of fixed positive
characteristic with divisible or $p$-divisible value groups and fixed
finite residue field.
%

\begin{theorem}                             \label{sptame}
a) \ Every elementary class {\bf C} of tame fields of fixed positive
characteristic with divisible value group and fixed residue field $\F_q$
(where $q=p^n$ for some prime $p$ and some $n\in\N$) is model complete,
complete and decidable. Moreover, it possesses a model of transcendence
degree 1 over $\F_q$ that admits an elementary embedding in every
$\aleph_1$-saturated member of {\bf C}.
\sn
b) \ If ``divisible value group'' is replaced by ``value group
elementarily equivalent to $\frac{1}{p^{\infty}}\Z$'', then {\bf C}
remains elementary, complete and decidable.
%
%
\end{theorem}
\begin{proof}
a): \ The theory of divisible ordered abelian groups is model complete,
complete and decidable, cf.\ \cite{[Ro--Zk]} (note that model
completeness and decidability are not explicitly stated in the theorems,
but follow from their proofs). The same holds trivially for the theory
of the finite field $\F_q$ which has only $\F_q$ as a model (up to
isomorphism). Hence, model completeness, completeness and decidability
follow readily from Theorem~\ref{tamemod1} and Theorem~\ref{dec}. The
prime model is constructed as follows: The valued field $(\F_q(t),v_t)$
has value group $\Z$ and residue field $\F_q\,$. By adjoining suitable
roots of $t$ we can build an algebraic extension $(F',v_t)$ with value
group $\Q$ and residue field $\F_q$. Now we let $(F,v_t)$ be a maximal
immediate algebraic extension of $(F',v_t)$. By Theorem~\ref{tame}, it is
a tame field. Moreover, it admits $\{t\}$ as a standard valuation
transcendence basis over its prime field. Note that $|F| = \aleph_0$.
Since $\Q$ is a prime model of the theory of nontrivial divisible
ordered abelian groups, Embedding Lemma II (Lemma~\ref{EII}) shows that
$(F,v_t)$ admits an embedding in every $\aleph_1$-saturated member of
{\bf C}. By the model completeness that we have already proved, this
embedding is elementary.
\sn
b): \ The theory of $\frac{1}{p^{\infty}}\Z$ is clearly complete, and it
is decidable (and {\bf C} is still elementary) because it can be
axiomatized by a recursive set of elementary axioms.
Now the proof proceeds as in part a), except that we replace $\Q$ by
$\frac{1}{p^{\infty}}\Z$ and note that the latter admits an elementary
embedding in every elementarily equivalent ordered abelian group (again,
cf.\ \cite{[Ro--Zk]}).
\end{proof}
Note that in the case of b), model completeness can be reinstated by
adjoining a constant symbol to the language and adding axioms that state
that the value of the element named by this symbol is divisible by no
prime but $p$.

%
%
\subsection{Separably defectless and separably tame fields}
\label{sectstsd}
We prove part a) of Theorem~\ref{septameAKE}:
%
\sn
Assume that $vK\ec vL$ and $Kv\ec Lv$. Since $vK$ is cofinal in $vL$, we
know that $(K,v)^c$ is contained in $(L,v)^c$. The compositum
$(L.K^c,v)$, taken in the completion $(L,v)^c$, is an immediate
extension of $(L,v)$. Thus, $vK^c=vK\ec vL=vL.K^c$ and
$K^cv=Kv\ec Lv=(L.K^c)v$.
%
%
Since $(K,v)$ is a henselian separably defectless field,
$(K,v)^c$ is henselian by Theorem~32.19 of \cite{[W]} and defectless by
Theorem~5.2 of \cite{[K8]}. As $(L|K,v)$ is an extension without
transcendence defect, the same holds for $(L.K^c|K^c,v)$; indeed, every
subextension of $L.K^c|K^c$ of finite transcendence degree is contained
in $L'.K^c|K^c$ for some subextension $L'|K$ of finite transcendence
degree, and since $(K^c|K,v)$ is immediate, a standard valuation
transcendence basis of $(L'|K,v)$ is also a standard valuation
transcendence basis of $(L'.K^c|K^c,v)$. By Theorem~\ref{wtA}, it now
follows that
\[
(K^c,v)\ec (L.K^c,v)\;.
\]
By Proposition~\ref{ecdec}, this implies that $(K,v)\ec (L,v)$.
\qed
\sn
%
%
%

\parm
We can now prove part b) of Theorem~\ref{septameAKE}:
%
\sn
Assume that $(K,v)$ is separably tame and that $(L|K,v)$ is a separable
extension with $vK\ec vL$ and $Kv\ec Lv$. If $\chara K=0$, then $(K,v)$
is tame and we have already proved that $(L|K,v)$ satisfies the
AKE$^\exists$ Principle. So we assume that $\chara K=p>0$. The perfect
hull $K^{1/p^{\infty}}$ of $K$ admits a unique extension $v$ of the
valuation of $K$, and with this valuation it is a subextension of the
completion of $K$, according to Lemma~\ref{Xsept}. In particular,
$(K^{1/p^{\infty}}|K,v)$ is immediate.
By Lemma~\ref{XKsK}, $(K^{1/p^{\infty}},v)$ is a tame field. Both
$K^{1/p^{\infty}}$ and $L.K^{1/p^{\infty}}$ are subfields of the perfect
hull $(L^{1/p^{\infty}},v)$ of $(L,v)$, whose value group is the
$p$-divisible hull of $vL$ and whose residue field is the perfect hull
of $Lv$.
%
%
As $vK=vK^{1/p^{\infty}}$ is $p$-divisible and $Kv=Kv^{1/p^{\infty}}$
is perfect, Lemma~\ref{te} shows that our side conditions yield that
$vK^{1/p^{\infty}} \ec v(L.K^{1/p^{\infty}})$ and $Kv^{1/p^{\infty}}\ec
(L.K^{1/p^{\infty}})v$. According to the AKE$^\exists$ Principle for
tame fields (Theorem~\ref{tameAKE}), this yields that
\[
(K^{1/p^{\infty}},v)\ec (L.K^{1/p^{\infty}},v)\;.
\]
By Proposition~\ref{ecdec}, this implies that $(K,v)\ec (L,v)$ since
$K$ is dense in $(K^{1/p^{\infty}},v)$.
\qed

\mn
Related to these results are results of \ind{F.~Delon} \cite{[D]}. She
showed that the \bfind{elementary class of algebraically maximal
Kaplansky fields of fixed $p$-degree} is relatively complete. Adding
predicates to the language of valued fields which guarantee that every
extension is separable, she also obtained relative model completeness.
We will discuss the case of separably tame fields of fixed $p$-degree in
a subsequent paper.

%
%
%
\newcommand{\lit}[1]{\bibitem{#1}}

\end{document}